\documentstyle[12pt,times,graphicx,subfigure,pstricks,amssymb]{article}
\newcommand\beq[1]{ \begin{equation}\label{#1} }
\newcommand{\eeq}{ \end{equation} }
\newcommand{\beqno}{ \[ }
\newcommand{\eeqno}{ \] }
\newcommand\beqa[1]{ \begin{eqnarray} \label{#1}}
\newcommand{\eeqa}{ \end{eqnarray} }
\newcommand{\beqano}{ \begin{eqnarray*} }
\newcommand{\eeqano}{ \end{eqnarray*} }
\newcommand\equ[1]{{\rm (\ref{#1})}}
\newcommand{\eps}{\varepsilon}
\newcommand{\R}{\mathbb{R}}
\newcommand{\Z}{\mathbb{Z}}

\newcommand{\T}{\mathbb{T}}
\newcommand{\C}{\mathbb{C}}

\begin{document}

\title{Stability estimates of nearly--integrable systems with dissipation and non--resonant frequency}

\author{ Alessandra Celletti\\
{\footnotesize Dipartimento di Matematica}\\
{\footnotesize Universit\`a di Roma Tor Vergata}\\
{\footnotesize Via della Ricerca Scientifica 1}\\
{\footnotesize I-00133 Roma (Italy)}\\
{\footnotesize \texttt{(celletti@mat.uniroma2.it)}}
\and
Christoph Lhotka\\
{\footnotesize Dipartimento di Matematica}\\
{\footnotesize Universit\`a di Roma Tor Vergata}\\
{\footnotesize Via della Ricerca Scientifica 1}\\
{\footnotesize I-00133 Roma (Italy)}\\
{\footnotesize \texttt{(lhotka@mat.uniroma2.it)}}
}

\maketitle

$\qquad\qquad\qquad\qquad\qquad\qquad\qquad\qquad$ \sl To Professor Anatoly Pavlovich Markeev,

$\qquad\qquad\qquad\qquad\qquad\qquad\qquad\qquad$ on the occasion of his 70th birthday \rm

\vglue1cm

\begin{abstract}
We consider a dissipative vector field which is represented by a
nearly--integrable Hamiltonian flow to which a non symplectic
force is added, so that the phase space volume is not preserved. The
vector field depends upon two parameters, namely the perturbing
and dissipative parameters, and by a drift function. We study the
general case of an $\ell$--dimensional, time--dependent vector
field. Assuming to start with non--resonant initial conditions, we
prove the stability of the variables which are actions of the conservative system
(namely, when the dissipative parameter is set to zero)
for exponentially long times. In order to construct the normal form,
a suitable choice of the drift function must be performed. We also provide
some simple examples in which we construct explicitly the normal form, we make a comparison
with a numerical integration and we compute theoretical bounds on the parameters
as well as we give explicit stability estimates.
\end{abstract}

\vglue1cm

\noindent \bf Keywords. \rm Dissipative system, Stability, Non--resonant motion.

\vglue1cm

\tableofcontents

\section{Introduction}\label{sec:intro}
A breakthrough in the theory of the stability of
nearly--integrable Hamiltonian systems was achieved by the seminal
works of A.N. Kolmogorov (\cite{K}), V.I. Arnold (\cite{Ar1},
\cite{Arenc}), J. Moser (\cite{M}) and N.N. Nekhoroshev
(\cite{nekh}, \cite{nekh2}). Nearly--integrable systems can be
modeled by Hamiltonian functions of the form \beq{ham}
H(y,x,t)=H_0(y)+\varepsilon H_1(y,x,t)\ , \eeq where, for an
$\ell$--dimensional system, $y\in{\R}^\ell$,
$(x,t)\in{\T}^{\ell+1}$, $H_0$, $H_1$ are regular functions, $H_1$
being periodic in $x$ and $t$, and $\varepsilon$ is a small parameter,
called the \sl perturbing parameter. \rm As far as
$\varepsilon=0$, Hamilton's equations associated to \equ{ham} are
integrable. The actions are constants and the motion is described
by periodic or quasi--periodic solutions, according to the
rational or irrational character of the frequency vector, say
$\omega(y)=\partial H_0(y)/\partial y$. In this context, under
general assumptions, KAM theory yields the persistence of
invariant tori for sufficiently small values of the perturbing
parameter, which implies a stability property for low--dimensional
systems, while higher dimensional models may admit diffusion through
the invariant tori. A stability result for arbitrary dimensions
can be obtained through Nekhoroshev's theorem, which guarantees
the stability of the actions (namely the confinement in a given
domain of the action space) for exponentially long times. Both KAM
and Nekhoroshev's theorems are constructive, in the sense that
they can be explicitly implemented to provide bounds on the
parameters ensuring the persistence of invariant tori or the
confinement of the actions over exponential times. For this reason
the theorems have been widely used to investigate physical systems
in different contexts. Most notably, Celestial Mechanics was a
spur for the development of analytical results about the stability
of nearly--integrable systems. In this field, there exist many
applications of KAM theorem (compare with \cite{CE10} and references therein)
and of Nekhoroshev's theorem  (see, e.g., \cite{BFG}, \cite{BGF},
\cite{BiaC}, \cite{CF}, \cite{CG}, \cite{ES} \cite{GDFGS},
\cite{GS}, \cite{LED}); most of the latter papers deal with the
stability of the triangular Lagrangian points. All these results
are obtained in a conservative framework.
However, as it is well known, many interesting physical systems of
Celestial Mechanics are affected by a small dissipation. We quote,
for example, the three--body problem with Poynting--Robertson drag
(\cite{BF}, \cite{CSLF}) or the spin--orbit model (compare with \cite{M1}) with tidal
torque (\cite{CE10}). These problems can be modeled by a
nearly--integrable system to which a small dissipation is added
(see, e.g., \cite{BST}).\\

\noindent While in the conservative setting we find periodic
orbits, invariant tori as well as chaotic motions, in the
dissipative context (where the
phase space volume is not preserved by time evolution) we speak of periodic attractors,
quasi--periodic attractors and strange attractors.
Beside the case of contracting and expanding systems, we consider
also \sl oscillating \rm systems for which the energy varies periodically
around a mean value. Coexisting attractors have been
established in, e.g., \cite{FGHY}, \cite{CE10}; periodic
attractors are shown to exist within parameter regions known as
Arnold's tongues (\cite{BST}, see also \cite{CD}). Concerning
invariant attractors, KAM results (both analytical and numerical)
provide their persistence under general assumptions (see \cite{broer}, \cite{CC},
\cite{CCL11}, \cite{CCdiss}); an application of the converse KAM
theorem (about the non--existence of quasi--periodic
attractors) has been presented in \cite{CMK}.\\

\noindent In this paper we exploit exponential type estimates for
nearly--integrable dissipative vector fields, which are defined as follows: we
consider a vector field depending on two parameters, say
$\varepsilon$ and $\mu$, such that for $\mu=0$ we obtain a
Hamiltonian vector field in action--angle variables, say
$(y,x)\in{\R}^\ell\times{\T}^\ell$; we consider the most general
case of a vector field depending also periodically on the time $t$. We
assume that for $\mu\not=0$ the vector field is not conservative;
at the first order in $\mu$ the vector field contains a drift
function, say $\eta=\eta(y,x,t)$, which must be suitably chosen in
order to meet some compatibility conditions, which allow to construct
a proper normal form. In this context we
establish stability estimates by proving the confinement
of the actions of the integrable approximation within a given
domain of the phase space for exponentially long times (see
\cite{benetal}, \cite{gionekh} and especially \cite{poe1}
which is at the basis of the present work) in the non--conservative
non--integrable system. The proof relies on the construction of a
proper normal form, which results from the composition of a
conservative transformation (removing the conservative terms to
suitable orders in $\varepsilon$) and a dissipative transformation,
which acts on the $\mu$--dependent terms. Perturbative methods for
vector fields have also been developed in \cite{fasso} (see also \cite{DGL},
\cite{LI05}, \cite{LI10}). The case of resonant
initial conditions is treated in \cite{CLres}, where one needs to construct a Lyapunov function (i.e. the
energy function associated to the conservative system with $\mu=0$), which must
be used in order to bound  the variation of the normal form coordinates;
in order to achieve the result, the resonant case requires to work in the
extended phase space and to impose  the quasi--convexity assumption, which
is not needed in the non--resonant case treated here. The final results are conceptually
different, since in the non--resonance case we obtain a stability result valid
for exponential times, while in the resonant case of \cite{CLres} the stability time may vary
linearly or exponentially with the parameters, according to the structure of the
vector field and to the choice of the resonance condition.

\noindent
We believe that an explicit construction of the normal form and of the stability
result in the non--resonant case is certainly of great interest in view of applications to
concrete models.
Just to quote a field which is familiar to the authors, a  non--resonant application in
Celestial Mechanics can be of interest in order to bound the motion of the majority of the
main belt asteroids.

\noindent
In this paper we provide also concrete
examples by investigating two model problems: a strictly dissipative vector field and a system with
oscillating energy. We test the accuracy of the normal form
construction by comparing the analytical expression with a
numerical integration (see \cite{CLInt} for a discussions about
the integration algorithms in nearly--Hamiltonian systems). Furthermore, we
compute concrete estimates on the parameters (namely, the radius
of the action domain and the stability time) and we investigate
their dependence on quantities like the normal form order, obtaining
stability times which grow exponentially with the normalization order.

\noindent
This paper is organized as follows. In Section~\ref{sec:setup} we
set--up the definitions and notations. In
Section~\ref{sec:stability} we state the normal form Lemma and the
main Theorem; the proofs are given in Section~\ref{sec:proof}. An
application of the normal form Lemma and the computation of the
stability estimates according to the results of the Theorem are
given in Section~\ref{sec:example} .

\section{Set--up: notations and assumptions}\label{sec:setup}
We consider the following $\ell$--dimensional, time--dependent vector field
\beqa{SE}
\dot x&=&\omega(y)+\varepsilon h_{10,y}(y,x,t)+\mu f_{01}(y,x,t)\nonumber\\
\dot y&=&-\varepsilon h_{10,x}(y,x,t)+\mu
\Big(g_{01}(y,x,t)-\eta(y,x,t)\Big)\ , \eeqa where $y\in A$ with
$A$ being an open domain of $\R^\ell$, $(x,t)\in{\T}^{\ell+1}$,
$\varepsilon\in{\R}_+$, $\mu\in{\R}\backslash\{0\}$, $\omega$ and $\eta$ are
real--analytic, $\ell$--dimensional vector functions with
components $(\omega^{(1)},...,\omega^{(\ell)})$ and
$(\eta^{(1)},...,\eta^{(\ell)})$. We assume that $h_{10}$,
$f_{01}$, $g_{01}$ are known periodic, real--analytic,
$\ell$--dimensional functions defined on $A\times {\T}^{\ell+1}$.
In all this paper we adopt the following notations and definitions.

$(i)$ The subscripts $x$, $y$ denote derivatives with respect to $x$, $y$.

$(ii)$ For integers $j$, $m$, the symbol $F_{jm}$ denotes that the function $F_{jm}$ is of order
$\varepsilon^j\mu^m$.

$(iii)$ We say that a function $F$ is of order $k$ in
$\varepsilon$ and $\mu$, in symbols $F\in O_k(\varepsilon,\mu)$, if its
Taylor series expansion in $\varepsilon$, $\mu$ contains only
powers of $\varepsilon^j\mu^m$ with $j+m\geq k$.

$(iv)$ For any integer vector
$k=(k_1,...,k_\ell)\in{\Z}^\ell$ we introduce the norm $|k|\equiv
|k_1|+...+|k_\ell|$.

$(v)$ We denote by a bar the average of a function over the angle
variables, while the tilde denotes the oscillatory part; more specifically,
we decompose a function $F=F(y,x,t)$ as
$$
F(y,x,t)=\bar F(y)+\tilde F(y,x,t)\ ,
$$
where the average $\bar F(y)$ is given by
$$
\bar F(y)\equiv \langle F(y,x,t)\rangle_{x,t}=
{1\over {(2\pi)^{\ell+1}}}\int_{{\T}^{\ell+1}} F(y,x,t)\ dx dt\ ,
$$
while the oscillatory part is defined as $\tilde F(y,x,t)\equiv F(y,x,t)-\bar F(y)$.

\vskip.1in

\noindent
We assume that there exists a
subset $D\subseteq A$ such that the vector function
$\omega=\omega(y)$ satisfies the following \sl non--resonance condition \rm
up to a suitable order $K$ with $K\in{\Z}_+$: \beq{NR}
|\omega(y)\cdot k+m|\geq a\quad {\rm for\ all\ } y\in D\ ,\qquad
k\in{\Z}^\ell\ ,\quad m\in{\Z}\ ,\quad |k|+|m|\leq K\ , \eeq where
$a$ is a strictly positive real constant and the dot denotes the scalar product.

\noindent
With reference to \equ{SE} we call $\varepsilon$ the \sl perturbing
parameter, \rm while we refer to $\mu$ as the \sl dissipative
parameter. \rm For $\varepsilon=\mu=0$ the equations \equ{SE} are trivially
integrated as
\beqano
x(t)&=&x(0)+\omega(y(0))t\nonumber\\
y(t)&=&y(0)\ ,
\eeqano
where $x(0),y(0)$ are the initial conditions at time $t=0$. For $\varepsilon\not=0$
small and $\mu=0$, the equations \equ{SE} reduce
to the conservative vector field associated to the
nearly--integrable Hamiltonian function
$$
H(y,x,t)=h_{00}(y)+\varepsilon h_{10}(y,x,t)\ ,
$$
where $h_{00}$ is such that $\omega(y)=h_{00,y}(y)$.  For $\mu\not=0$ we assume that the vector field is dissipative for any
$y\in A$, $(x,t)\in{\T}^{\ell+1}$, namely there exists a subset of
the phase space, which is contracted or expanded asymptotically by
time evolution into a compact set (\cite{ruelle}). We
also consider systems such that the energy is oscillating around a
mean value.

\noindent
We refer to $\eta=\eta(y,x,t)$ as the \sl drift \rm vector function with components
$(\eta^{(1)}(y,x,t),...,\eta^{(\ell)}(y,x,t))$, which can be expanded as
$$
\eta^{(k)}(y,x,t)=\sum_{m=0}^\infty \sum_{j=0}^m
\eta_{j,m-j+1}^{(k)}(y,x,t) \varepsilon^j\mu^{m-j}\ ,\qquad
k=1,...,\ell\ .
$$
We remark that $\eta$ is an unknown function, which will be properly chosen
in order to meet some compatibility requirements in order to perform a suitable normal form
(compare with KAM results like in \cite{CCL11}).

\noindent
Let a function $f=f(y,x,t)$ be defined for $y\in A$,
$(x,t)\in{\T}^{\ell+1}$; we denote by $C_{r_0}(A)$ the complex neighborhood of $A$
of radius $r_0$, namely
$$
C_{r_0}(A)\equiv \{{y\in\C^\ell:}\ \|y-y_A\|\leq r_0\ {\rm for\ all\ } y_A\in A\}\ ,
$$
where $\|\cdot\|$ denotes the Euclidean norm. Let $C_{s_0}({\T}^{\ell+1})$ be
the complex strip of radius $s_0$ around ${\T}^{\ell+1}$, namely
$$
C_{s_0}({\T}^{\ell+1})\equiv \{(x,t)\in{\C}^{\ell+1}:\ \max_{1\leq j\leq \ell}|\Im(x_j)|
\leq s_0\ ,\ |\Im(t)|\leq s_0\}\ ,
$$
where $\Im$ denotes the imaginary part. Let us denote the Fourier expansion of a function $f=f(y,x,t)$ as
$$
f(y,x,t)=\sum_{(k,m)\in{\Z}^{\ell+1}} \hat f_{km}(y) e^{i(k\cdot x+mt)}\ .
$$
For an analytic function on $A\times {\T}^{\ell+1}$ we introduce the norm
$$
\|f\|_{r_0,s_0}\equiv \sup_{y\in C_{r_0}(A)}\ \sum_{(k,m)\in{\Z}^{\ell+1}}
|\hat f_{km}(y)|\, e^{(|k|+|m|)s_0}\ ,
$$
while for a function $g=g(y)$, $g:{\R}^\ell\to \R$, we define $\|g\|_{r_0}\equiv\sup_{y\in C_{r_0}(A)} \|g(y)\|$.
For an $\ell$--dimensional vector function $v$ with components $(v_1,...,v_\ell)$, we define
$$
\|v\|_{r_0,s_0}\equiv \sqrt{\sum_{j=1}^\ell \|v_j\|_{r_0,s_0}^2}\ .
$$
For a function $f=f(y,x,t)$, $f:{\R}^\ell\times{\T}^{\ell+1}\to {\R}$, and for any positive integer $K$, we denote by $f^{\leq K}$, $f^{>K}$ the sum over the components with Fourier modes less or equal, or
respectively greater than $K$, namely
$$
f^{\leq K}(y,x,t)\equiv \sum_{(k,m)\in{\Z}^{\ell+1},\ |k|+|m|\leq K}
\hat f_{km}(y)e^{i(k\cdot x+mt)}
$$
and
$$
f^{> K}(y,x,t)\equiv \sum_{(k,m)\in{\Z}^{\ell+1},\ |k|+|m|> K}
\hat f_{km}(y)e^{i(k\cdot x+mt)}\ .
$$

\section{Stability for exponential times: statement of the results}\label{sec:stability}
Stability estimates for exponential times are obtained by implementing a change of
coordinates such that the vector field \equ{SE} is transformed to a suitable normal form of order $N$.
Precisely, let us consider a coordinate transformation close to the identity, say
$\Xi^{(n)}:A \times {\T}^{\ell+1} \to {\R}^\ell \times {\T}^{\ell+1}$, such that
\beq{c} (Y,X,t)=\Xi^{(N)}(y,x,t)\ ,\qquad Y\in{\R}^\ell\ ,\qquad (X,t)\in{\T}^{\ell+1}\ ,
\eeq where $\Xi^{(N)}$ depends parametrically also on $\varepsilon$, $\mu$;
we require that the transformation acts as the identity on the time.
Let us assume that on a suitable parameter domain, the coordinate transformation \equ{c} can be inverted as
$$
y=y(Y,X,t),\qquad x=x(Y,X,t)\ .
$$
In the following \sl Normal Form Lemma \rm we look for a change of
coordinates \equ{c} such that \equ{SE} is transformed to the following normal
form of order $N$:
\beqa{NF}
\dot X&=&\Omega_d^{(N)}(Y;\varepsilon,\mu)+F_{N+1}(Y,X,t)+F_{01}^{>K}(Y,X,t)\nonumber\\
&+&\varepsilon h^{>K}_{10,y}(y(Y,X,t),x(Y,X,t),t)+\mu f^{>K}_{01}(y(Y,X,t),x(Y,X,t),t)\nonumber\\
\dot Y&=&G_{N+1}(Y,X,t)+G_{01}^{>K}(Y,X,t)\nonumber\\
&-&\varepsilon
h^{>K}_{10,x}(y(Y,X,t),x(Y,X,t),t)+ \mu
g^{>K}_{01}(y(Y,X,t),x(Y,X,t),t)\ , \eeqa where
$F_{N+1}$, $G_{N+1}$ are vector functions of order
$O_{N+1}(\varepsilon,\mu)$;
$F_{01}^{>K}$, $G_{01}^{>K}$ are functions of order $O_1(\varepsilon,\mu)$,
depending on $\omega$, $h_{10}$, $f_{01}$, $g_{01}$ and on the normal
form transformation;
$\Omega_d^{(N)}:\R^\ell\rightarrow\R^\ell$ is an
$\ell$--dimensional vector function, related to $\omega(Y)$ by
$$
\Omega_d^{(N)}(Y;\varepsilon,\mu)\equiv
\omega(Y)+\sum_{m=0}^N\sum_{j=0}^m\Omega_{j,m-j}(Y)\varepsilon^j\mu^{m-j}\ ,
$$
being $\Omega_{j,m-j}(Y)$ known vector functions and being $\Omega_{00}=0$.
The meaning of \equ{NF} is that the normalized equations take the following form: the $X$--variation is provided by a modified
frequency $\Omega_d^{(N)}$ plus higher order terms in the normalization order $N$ or in the non--resonance order $K$;
the variation of the normal form variable $Y$ is constant, beside higher order terms in $N$ and $K$.

\noindent
The coordinate transformation $\Xi^{(N)}$ results from the composition of a transformation $\Xi^{(N)}_c$ acting on the conservative part and a change of coordinates $\Xi^{(N)}_d$ acting on the dissipative part, say:
\beq{7bis}
(Y,X,t)=\Xi^{(N)}_d \circ \Xi^{(N)}_c(y,x,t)\ .
\eeq
Setting the intermediate variables as $(\tilde y,\tilde x,t)\equiv \Xi^{(N)}_c(y,x,t)$, the
transformation $\Xi^{(N)}_c$ is implicitly defined through a
sequence of generating functions $\psi_{j0}=\psi_{j0}(\tilde
y,x,t)$, $j=1,...,N$, as \beqa{A}
\tilde x&=&x+\sum_{j=1}^N \psi_{j0,y}(\tilde y,x,t)\varepsilon^j\equiv x+\psi_y^{(N)}(\tilde y,x,t)\nonumber\\
y&=&\tilde y+\sum_{j=1}^N \psi_{j0,x}(\tilde
y,x,t)\varepsilon^j\equiv \tilde y+\psi_x^{(N)}(\tilde y,x,t)\ ,
\eeqa while $\Xi^{(N)}_d$ is defined by introducing suitable
functions $\alpha_{jm}$, $\beta_{jm}$, $j,m\in{\Z}_+$, as
\beqa{B} X&=&\tilde x+\sum_{m=0}^N\sum_{j=0}^m \alpha_{j,m-j}
(\tilde y,\tilde x,t)\varepsilon^j\mu^{m-j}\nonumber\\
Y&=&\tilde y+\sum_{m=0}^N\sum_{j=0}^m \beta_{j,m-j} (\tilde
y,\tilde x,t)\varepsilon^j\mu^{m-j}\ , \eeqa with the properties that
$\alpha_{i0}(\tilde y,\tilde x,t)=\beta_{i0}(\tilde y,\tilde
x,t)=0$ for any $i\geq 0$ and $\langle \alpha_{ij}(\tilde y,\tilde
x,t)\rangle_{\tilde x,t}=\langle \beta_{ij}(\tilde y,\tilde
x,t)\rangle_{\tilde x,t}=0$. The proof of the following lemma provides an algorithm to compute explicitly
the vector functions $\psi_{j0}$, $\alpha_{j m}$, $\beta_{j m}$, together with a suitable drift function $\eta$,
which allows to achieve the desired normal form.

\vskip.2in

\noindent \bf Normal Form Lemma. \sl Consider the real analytic vector
field \equ{SE} defined on $A\times{\T}^{\ell+1}$ with complex extension in
$C_{r_0}(A)\times C_{s_0}({\T}^{\ell+1})$ for some $r_0$, $s_0>0$.
Let $D\subseteq A$ be such that for any $y\in D$, the frequency $\omega=\omega(y)$ satisfies
\equ{NR} for some $K\in{\Z}_+$, $a>0$. Then, there exist
$\varepsilon_0$, $\mu_0>0$ depending on $r_0$, $s_0$, $K$ and the norms of $\omega$, $h$, $f$, $g$
and there exists $\eta=\eta(y,x,t)$, such that for $(\eps,|\mu|)\leq(\eps_0,\mu_0)$, one can find a coordinate
transformation close to the identity, say
$\Xi^{(N)}:A \times {\T}^{\ell+1} \to {\R}^\ell \times {\T}^{\ell+1}$
for a suitable normalization order $N\in\Z_+$, which brings \equ{SE} into a normal
form of order $N$ as in \equ{NF}. For $\tilde r_0<r_0$, the
normalized frequency $\Omega_d^{(N)}$ in \equ{NF}
is bounded by
$$
\|\Omega_d^{(N)}-\omega\|_{\tilde r_0}\leq C_\omega\lambda\ ,
$$
where $\lambda\equiv\max(\eps,|\mu|)$ and
$C_\omega$ is a suitable positive constant depending on $r_0$, $N$ and on the norms
of $\omega$, $h$, $f$, $g$, respectively. Let
$R_0<r_0$, $S_0<s_0$; with reference to the normal form equations
\equ{NF}, the following estimate holds: \beq{FG}
\|G_{N+1}\|_{R_0,S_0}+\|G^{>K}_{01}\|_{R_0,S_0}+\varepsilon \|h^{>K}_{10,x}\|_{R_0,S_0}+|\mu|
\|g^{>K}_{01}\|_{R_0,S_0} \leq
C_Y\lambda^{N+1}\ , \eeq for some positive constant $C_Y$ depending on $r_0$, $s_0$, $N$, $K$ and
the norms of $\omega$, $h$, $f$, $g$.
Choosing\footnote{The choice of $N$ can be performed as follows.
The relation $\lambda^N=e^{-K\tau_0}$ implies $N\log
\lambda=-K\tau_0$, namely $N=[K\tau_0/|\log \lambda|]$, where
$[\cdot]$ denotes the integer part.} $N=[K\tau_0/|\log \lambda|]$
for some $\tau_0>0$, one gets that \equ{FG} becomes
$$
\|G_{N+1}\|_{R_0,S_0}+\|G^{>K}_{01}\|_{R_0,S_0}+\varepsilon \|h^{>K}_{10,x}\|_{R_0,S_0}+|\mu|
\|g^{>K}_{01}\|_{R_0,S_0} \leq C_Y\lambda
e^{-K\tau_0}\ .
$$
Finally, denoting by $\Pi_y$ the projection on the
$y$--coordinate, one has
\beq{Pi} \|\Pi_y(\Xi_d^{(N)}\circ\Xi_c^{(N)})-Id\|\leq C_p\lambda\ ,
\eeq for some constant $C_p>0$ depending on $r_0$, $s_0$, $N$ and on the norm of
$\omega$, $h$, $f$ and $g$. \rm \\

\noindent
\bf Remark. \rm As an outcome of the proof of the Normal Form Lemma, the drift $\eta$ will depend only on
the normal form variable $Y$, i.e. $\eta=\eta(Y)$. The explicit form of $\eta$ depends upon
the functions $h_{10}$, $f_{01}$, $g_{01}$ appearing in \equ{SE}; we remark that the value of $\eta$
is determined by the requirement that $Y$ is constant up to the normalization order $N$, say
$Y=Y_0+O_{N+1}(\varepsilon,\mu)$; on the other hand, $Y_0$ depends on $y_0$, which is chosen
so that the frequency $\omega=\omega(y)$ satisfies \equ{NR}. The fact that $\eta$ is linked to
the form of the vector field and to the frequency of motion appears also in KAM proofs for
dissipative (or conformally symplectic) systems (compare with \cite{CCdiss}, \cite{CCL11}).

\vskip.2in

\noindent Before giving the proof of the Lemma, we provide the
statement of the main result, namely the confinement of the $y$
variables for exponential times, which will be obtained through the
Normal Form Lemma under the non--resonance condition \equ{NR}.

\vskip.2in

\noindent
\bf Theorem. \sl Consider the vector field \equ{SE} defined on $A\times {\T}^{\ell+1}$,
and let $D\subseteq A$ be such that for any $y\in D$ the frequency
$\omega=\omega(y)$ satisfies \equ{NR}.
Assume there exists $\varepsilon_0$, $\mu_0$ such that for
$(\eps,|\mu|)\leq (\eps_0,\mu_0)$ the Normal Form Lemma holds.
Then, there exist positive parameters $\rho_0,\tau_0>0$, such that for every
solution $(y(t),x(t))$ at time $t>0$ with initial position $(y(0),x(0))\in D\times {\T}^{\ell}$ one
has for $\lambda\equiv\max(\eps,|\mu|)$:
$$
\|y(t)-y(0)\|\leq \rho_0\lambda \quad {for}\quad t\leq C_te^{K\tau_0}\ ,
$$
for some positive constant $C_t$, where $\rho_0$ and $C_t$ depend on $r_0$, $s_0$, $N$, $K$,
$\eps_0$, $\mu_0$ and on the norms of $\omega$, $h$, $f$, $g$. \rm

\vskip.1in

\noindent
\bf Remark. \rm Notice that we obtain the standard formulation of the stability time in terms of
an exponential estimate in the inverse of the small parameters, by adopting a proper choice of
$K\tau_0$, say $K\tau_0\leq ({1\over \lambda})^c$ for some constant $c>0$;
in this case one has that the stability time estimate is $t\leq C_t\, e^{{1\over \lambda_0^c}\,({\lambda_0\over \lambda})^c}$
for $\lambda_0\equiv \max(\eps_0,\mu_0)$.

\section{Proof of the Normal Form Lemma and of the Theorem}\label{sec:proof}
In this section we first outline the general scheme of the proof and then we provide the
complete proof of the Normal Form Lemma, followed by that of the main Theorem. For easiness of readability
technical Lemmas appear in the Appendixes. We start by implementing a coordinate change of
variables of the form \equ{7bis}. In particular, we define the intermediate variables $(\tilde y,\tilde
x,t)\in\R^\ell\times$ ${{\T}^{\ell+1}}$, which provide the transformation
\equ{A} in order that the following conservative normal form is obtained:
\beqa{NCF}
\dot{\tilde x}&=&\Omega_c^{(N)}(\tilde y;\varepsilon)+\mu \sum_{j=1}^N F_{j1}^{(1,N)}(\tilde y,\tilde x,t)\varepsilon^j
+\sum_{j=1}^N F_{j0}^{(1,N,>K)}(\tilde y,\tilde x,t)\varepsilon^j\nonumber\\
&+&\varepsilon h^{>K}_{10,y}(y(\tilde y,\tilde x,t),x(\tilde y,\tilde x,t),t)+
\mu f^{>K}_{01}(y(\tilde y,\tilde x,t),x(\tilde y,\tilde x,t),t)
+O_{N+1}(\varepsilon,\mu)\nonumber\\
\dot{\tilde y}&=&\mu F^{(2,N)}(\tilde y,\tilde x,t)-\mu\left(\eta_{N-1,1}(\tilde y,\tilde x,t)
\varepsilon^{N-1}+...+\eta_{0N}(\tilde y,\tilde x,t)\mu^{N-1}\right)\nonumber\\
&-&\varepsilon h^{>K}_{10,x}(y(\tilde y,\tilde x,t),x(\tilde y,\tilde x,t),t)+
\mu g^{>K}_{01}(y(\tilde y,\tilde x,t),x(\tilde y,\tilde x,t),t)\nonumber\\
&+&\sum_{j=1}^N F_{j0}^{(3,N,>K)}(\tilde y,\tilde x,t)\varepsilon^j+O_{N+1}(\varepsilon,\mu)\ ,
\eeqa
for suitable functions $F_{j1}^{(1,N)}$, $F_{j0}^{(1,N,>K)}$, $F^{(2,N)}$, $F_{j0}^{(3,N,>K)}$,
and being $\Omega_c^{(N)}(\tilde y;\varepsilon)\equiv
\omega(\tilde y)+\sum_{j=1}^N \Omega_{j0}(\tilde y)\varepsilon^j$,
where $\Omega_{j0}(\tilde y)$ can be explicitly determined. We
denote the inversion of \equ{A} as
\beqa{inv}
x&=&x(\tilde y,\tilde x,t)=\tilde x+\Gamma^{(x,N)}(\tilde y,\tilde x,t)\nonumber\\
y&=&y(\tilde y,\tilde x,t)=\tilde y+\Gamma^{(y,N)}(\tilde y,\tilde x,t)\ ,
\eeqa
which provides
the transformation $\Xi_c^{(N)}$. We will see that the generating function
$\psi_{j0}(\tilde y,\tilde x,t)$ in \equ{A} at the generic order $j$ is
the solution of an equation of the following form, defined in
terms of the intermediate set of variables:
$$
\omega(\tilde y)\ \psi_{j0,x}(\tilde y,\tilde x,t)
+\psi_{j0,t}(\tilde y,\tilde x,t)+\tilde L^{(j,\leq K)}(\tilde y,\tilde x,t)=0\ ,
$$
for a suitable known function $\tilde L^{(j,\leq K)}(\tilde y,\tilde x,t)$ with zero average over $(\tilde x,t)$.
The above equation can be solved provided $\omega=\omega(\tilde y)$ satisfies a
non--resonance condition of the form
$$
k\cdot \omega(\tilde y)+m\not=0\qquad {\rm for\ all}\ k\in{\Z}^\ell\ ,\ \ m\in{\Z}\ ,\
|k|+|m|\leq K\ ,
$$
which is guaranteed by \equ{NR} on a suitable domain.

\noindent After the implementation of the conservative
transformation, to achieve the normal form \equ{A} we
construct a change of coordinates $(Y,X,t)=\Xi_d^{(N)}(\tilde
y,\tilde x,t)$ defined as in \equ{B}, which allows to obtain the normal form \equ{NF}.
We will see that the functions $\beta_{j m}$ must satisfy an equation of
the form
\beq{beta}
\omega(Y) \beta_{j m,x}(Y,X,t)+\beta_{j m,t}(Y,X,t)+ N_{j m}^{\leq K}(Y,X,t)-\eta_{jm}(Y,X,t)=0\ ,
\eeq
for some known function $N_{j m}(Y,X,t)$; therefore, equation \equ{beta} can be
solved provided that the drift components $\eta_{jm}(Y)$ are chosen as the averages
$\bar N_{j m}$ of $N_{j m}^{\leq K}$:
$$
\eta_{jm}(Y,X,t)\equiv \eta_{jm}(Y)=\bar N_{j m}(Y)\ .
$$
The normal form equation \equ{beta} can be solved
provided that the frequency satisfies the non--resonance condition $k\cdot\omega(Y)+m\not=0$ for all
$k\in{\Z}^\ell$, $m\in{\Z}$, $|k|+|m|\leq K$, which is guaranteed by \equ{NR} on a suitable domain.
Once $\beta_{jm}$ is determined, we can proceed to compute $\alpha_{jm}$ by solving a normal form
equation again of the form \equ{beta}, but having zero average.

\vskip.1in

\noindent \bf Proof of the Normal Form Lemma. \rm We prove by
induction on the normal form order that we can determine the
transformations \equ{A} and \equ{B} so to obtain the normal form \equ{NF}. We start by
constructing the first order normal form through the implementation of the
conservative and then of the dissipative transformation; in a similar
way we construct the transformations for the order $N$. Being the
proof quite long, for sake of exposition we split it into separate steps.

\vskip.1in

\noindent
\bf Step 1: Conservative transformation for $N=1$. \rm

\noindent
Let us start with the conservative transformation for $N=1$, namely we implement the first order
change of variables
\beqa{deltac1}
\tilde x&=&x+ \varepsilon \psi_{10,y}(\tilde y,x,t)\nonumber\\
y&=&\tilde y+ \varepsilon \psi_{10,x}(\tilde y,x,t)\ ,
\eeqa
where $\psi_{10}=\psi_{10}(\tilde y,x,t)$ must be determined. Let $\tilde r_0<r_0$, $\delta_0<s_0$, $\tilde s_0\equiv s_0-\delta_0$; then \equ{deltac1} can be inverted as
\beqano
x&=&\tilde x+\varepsilon \Gamma^{(x,1)}(\tilde y,\tilde x,t)\nonumber\\
y&=&\tilde y+\varepsilon \Gamma^{(y,1)}(\tilde y,\tilde x,t)\ ,
\eeqano
for suitable functions $\Gamma^{(x,1)}$ and $\Gamma^{(y,1)}$, provided the following condition is satisfied
(see Appendix~A):
\beq{C1}
70\,\varepsilon \ \|\psi_{10,y}\|_{\tilde r_0,s_0}\ e^{2s_0}\delta_0^{-1}<1\ .
\eeq
Using \equ{deltac1} and \equ{SE}, we compute the time derivatives of $\tilde x$,
$\tilde y$ as
\beqa{star}
\dot{\tilde x}&=&\omega(\tilde y)+\varepsilon \omega_y(\tilde y)\psi_{10,x}(\tilde y,\tilde x,t)+
\varepsilon \tilde h_{10,y}(\tilde y,\tilde x,t)+\varepsilon\bar h_{10,y}(\tilde y)
+\mu f_{01}(\tilde y,\tilde x,t)\nonumber\\
&+&\varepsilon \omega(\tilde y)\psi_{10,yx}(\tilde y,\tilde x,t)+
\varepsilon \psi_{10,yt}(\tilde y,\tilde x,t)+O_2(\varepsilon,\mu)\nonumber\\
\dot{\tilde y}&=&-\varepsilon\tilde h_{10,x}(\tilde y,\tilde x,t)+\mu \Big(g_{01}(\tilde y,\tilde x,t)
-\eta_{01}(\tilde y,\tilde x,t)\Big)-\varepsilon\omega(\tilde y)\psi_{10,xx}(\tilde y,\tilde x,t)\nonumber\\
&-&\varepsilon\psi_{10,xt}(\tilde y,\tilde x,t)+O_2(\varepsilon,\mu)\ .
\eeqa
The conservative normal form is obtained by imposing that $\psi_{10}(\tilde y,\tilde x,t)$
satisfies the following normal form equations:
\beqa{otto}
\omega_y(\tilde y)\psi_{10,x}(\tilde y,\tilde x,t)+\psi_{10,yt}(\tilde y,\tilde x,t)+
\omega(\tilde y)\psi_{10,yx}(\tilde y,\tilde x,t)+\tilde h^{\leq K}_{10,y}(\tilde y,\tilde x,t)&=&0\nonumber\\
\omega(\tilde y)\psi_{10,xx}(\tilde y,\tilde x,t)+\psi_{10,xt}(\tilde y,\tilde x,t)+
\tilde h^{\leq K}_{10,x}(\tilde y,\tilde x,t)&=&0\ .
\eeqa
As a consequence, setting
\beq{omegag}
\Omega_c^{(1)}(\tilde y;\varepsilon)\equiv \omega(\tilde y)+\varepsilon \bar h_{10,y}(\tilde y) \ ,
\eeq
equations \equ{star} become
\beqa{punto}
\dot{\tilde x}&=&\Omega_c^{(1)}(\tilde y;\varepsilon)+\varepsilon h^{>K}_{10,y}(\tilde y,\tilde x,t)
+\mu f_{01}(\tilde y,\tilde x,t)+O_2(\varepsilon,\mu)\nonumber\\
\dot{\tilde y}&=&-\varepsilon h^{>K}_{10,x}(\tilde y,\tilde x,t)+
\mu \Big(g_{01}(\tilde y,\tilde x,t)-\eta_{01}(\tilde y,\tilde x,t)\Big)+O_2(\varepsilon,\mu)\ ,
\eeqa
which are recognized as being of the form \equ{NCF}.
We remark that equations \equ{otto} are obtained taking, respectively, the derivatives with respect to $y$
and $x$ of
$$
\omega(\tilde y)\psi_{10,x}(\tilde y,\tilde x,t)+\psi_{10,t}(\tilde y,\tilde x,t)+
\tilde h^{\leq K}_{10}(\tilde y,\tilde x,t)=0\ .
$$
Expanding $\tilde h^{\leq K}_{10}(\tilde y,\tilde x,t)$ in Fourier series as
$$
\tilde h^{\leq K}_{10}(\tilde y,\tilde x,t)= \sum_{(k,m)\in{\Z}^{\ell+1},|k|+|m|\leq K}
{\hat {\tilde h}_{10,km}(\tilde y)}e^{i(k\cdot\tilde x+mt)}\ ,
$$
where ${\hat {\tilde h}_{10,km}(\tilde y)}$ denote the Fourier coefficients,
the solution for $\psi_{10}(\tilde y,\tilde x,t)$ is given by:
$$
\psi_{10}(\tilde y,\tilde x,t)=i \sum_{(k,m)\in{\Z}^{\ell+1},|k|+|m|\leq K}
{{\hat {\tilde h}_{10,km}(\tilde y)}\over {\omega(\tilde y)\cdot k+m}}\
e^{i(k\cdot\tilde x+mt)}\ .
$$
To control the small divisors appearing in the previous expression, let us invert the
second in \equ{deltac1} as $\tilde y=y+\varepsilon R^{(1)}(y,x,t)$ for a suitable
function $R^{(1)}=R^{(1)}(y,x,t)$, provided that for $\tilde r_0'<r_0$ one has (compare with Appendix~A)
\beq{C2}
70\,\varepsilon\|\psi_{10,x}\|_{\tilde r_0,s_0}{1\over {\tilde r_0-\tilde r_0'}}<1\ .
\eeq
Then, we have that
$$
|\omega(\tilde y)\cdot k+m|\geq a-\varepsilon K\|R^{(1)}\|_{\tilde r_0',s_0} \|\omega_y\|_{r_0}\geq {a\over 2}\ ,
$$
if (compare with Appendix~A)
\beq{C3}
\varepsilon \leq {a\over {2K\|R^{(1)}\|_{\tilde r_0',s_0} \|\omega_y\|_{r_0}}}\ .
\eeq

\vskip.1in

\noindent
\bf Step 2: Dissipative transformation for $N=1$. \rm

\noindent
We proceed to reduce to normal form the dissipative part through a first--order
transformation of coordinates $\Delta_d^{(1)}$, which we write in components as
\beqa{D1}
X&=&\tilde x+\alpha_{01}(\tilde y,\tilde x,t)\mu\nonumber\\
Y&=&\tilde y+\beta_{01}(\tilde y,\tilde x,t)\mu\ ,
\eeqa
for some functions $\alpha_{01}$ and $\beta_{01}$ to be determined as follows. We premise that equations \equ{D1} can be inverted as
\beqano
\tilde x&=&X+\Delta^{(x,1)}(Y,X,t)\mu\nonumber\\
\tilde y&=&Y+\Delta^{(y,1)}(Y,X,t)\mu
\eeqano
for suitable functions $\Delta^{(x,1)}$ and $\Delta^{(y,1)}$ provided the following conditions are satisfied
(see Appendix~A):
\beqa{C4}
70\,|\mu|\ \|\alpha_{01}\|_{\tilde r_0,\tilde s_0}\ e^{2\tilde s_0}\tilde\delta_0^{-1}<1\nonumber\\
70\,|\mu|\ (\|\beta_{01}\|_{\tilde r_0,\tilde s_0}+|\mu|\|\beta_{01,x}\|_{\tilde r_0,\tilde s_0}\ \|\alpha_{01}\|_{\tilde r_0,\tilde s_0})\ {1\over {\tilde r_0-R_0}}<1\ ,
\eeqa
where $\tilde \delta_0\equiv \tilde s_0/2$, $R_0<\tilde r_0$ and
being $\|\Delta^{(x,1)}\|_{R_0,S_0}\leq\|\alpha_{01}\|_{\tilde r_0,\tilde s_0}$ with $S_0<\tilde s_0-\tilde \delta_0$.
Up to the second order, the inversion of \equ{D1} provides
\beqa{starstar}
\tilde x&=&X-\alpha_{01}(Y,X,t)\mu+O_2(\mu)\nonumber\\
\tilde y&=&Y-\beta_{01}(Y,X,t)\mu+O_2(\mu)\ .
\eeqa
Taking the derivative of \equ{starstar} and using \equ{punto}, we express $\dot X$, $\dot Y$ as a function of $X$, $Y$ as
\beqano
\dot X&=&\omega(Y)-\omega_y(Y)\beta_{01}(Y,X,t)\mu+\varepsilon{\bar h}_{10,y}(Y)\nonumber\\
&+&\varepsilon{h}^{>K}_{10,y}(Y,X,t)+\mu f_{01}(Y,X,t)+
\omega(Y)\alpha_{01,x}(Y,X,t)\mu\nonumber\\
&+&\alpha_{01,t}(Y,X,t)\mu+O_2(\varepsilon,\mu)\nonumber\\
\dot Y&=&-\varepsilon{h}^{>K}_{10,x}(Y,X,t)+
\mu \Big(g_{01}(Y,X,t)-\eta_{01}(Y,X,t)\Big)\nonumber\\
&+&\omega(Y)\beta_{01,x}(Y,X,t)\mu+\beta_{01,t}(Y,X,t)\mu
+O_2(\varepsilon,\mu)\ .
\eeqano
The dissipative normal form is obtained imposing that $\alpha_{01}$, $\beta_{01}$
and $\eta_{01}$ satisfy the following equations:
\beqa{XX}
\omega(Y)\alpha_{01,x}(Y,X,t)+\alpha_{01,t}(Y,X,t)-\omega_y(Y)\beta_{01}(Y,X,t)
+\tilde f_{01}^{\leq K}(Y,X,t)&=&0\nonumber\\
\omega(Y)\beta_{01,x}(Y,X,t)+\beta_{01,t}(Y,X,t)+\tilde g^{\leq K}_{01}(Y,X,t)+\bar g_{01}(Y)-\eta_{01}(Y)&=&0\ ,
\eeqa
where we have split $f_{01}$ into the sum of its average and of the oscillatory part,
namely
$$
f_{01}(Y,X,t)\equiv \bar f_{01}(Y)+\tilde f_{01}^{\le K}(Y,X,t)+f_{01}^{>K}(Y,X,t)
$$
and similarly for $g_{01}$: $g_{01}(Y,X,t)\equiv\bar g_{01}(Y)+\tilde g_{01}^{\le K}(Y,X,t)
+g_{01}^{>K}(Y,X,t)$.
From the first of \equ{XX} we see that the average of $\beta_{01}$ is zero and we can assume that also the
average of $\alpha_{01}$ is zero. Then, equations \equ{XX} can be solved,
provided that in the second equation the term $\eta_{01}(Y,X,t)$ is chosen so that
$$
\eta_{01}(Y,X,t)\equiv \eta_{01}(Y)=\bar g_{01}(Y)\ .
$$
The final normal form can be written as
\beqano
\dot X&=&\omega(Y)+\varepsilon{\bar h}_{10,y}(Y)+\mu \bar f_{01}(Y)
+\varepsilon h^{>K}_{10,y}(Y,X,t)+\mu f^{>K}_{01}(Y,X,t)+F_2(Y,X,t)\nonumber\\
\dot Y&=&-\varepsilon h^{>K}_{10,x}(Y,X,t)+\mu g^{>K}_{01}(Y,X,t)+G_2(Y,X,t)\ ,
\eeqano
where $F_2$, $G_2$ are $O_2(\varepsilon,\mu)$. These equations are recognized to be of the form
\equ{NF} with $\Omega_d^{(1)}(Y)=\omega(Y)+\varepsilon \bar h_{10}(Y)+\mu \bar f_{01}(Y)$.

\vskip.1in

\noindent
The solutions of \equ{XX} involves small divisors, which can be bounded as follows:
\beqno
|\omega(Y)\cdot k+m|\geq{a\over2}-|\mu| K\|\beta_{01}\|_{\tilde r_0,\tilde s_0}\|\omega_y\|_{r_0}>{a\over 4}\ ,
\eeqno
provided that (see Appendix~A)
\beq{C5}
|\mu| <{a\over {4K\|\beta_{01}\|_{\tilde r_0,\tilde s_0} \|\omega_y\|_{r_0}}}\ .
\eeq

\vskip.1in

\noindent
\bf Step 3: Conservative transformation for the order $N$. \rm

\noindent
Assuming that the Lemma holds up to the order $N-1$, we prove it for the order $N$, starting from the change of
variables \equ{A}, that we invert as
\beqa{deltac3}
x&=&\tilde x+\sum_{j=1}^{N} \Gamma_{j0}^{(x)}(\tilde y,\tilde x,t)\varepsilon^j-
\psi_{N0,y}(\tilde y,\tilde x,t)\varepsilon^N+O_{N+1}(\varepsilon)
\equiv \tilde x+\Gamma^{(x,N)}(\tilde y,\tilde x,t)\nonumber\\
y&=&\tilde y+\sum_{j=1}^{N} \Gamma_{j0}^{(y)}(\tilde y,\tilde x,t)\varepsilon^j+
\psi_{N0,x}(\tilde y,\tilde x,t)\varepsilon^N+O_{N+1}(\varepsilon)\equiv \tilde y+
\Gamma^{(y,N)}(\tilde y,\tilde x,t)\ ,
\eeqa
where $\Gamma_{j0}^{(x)}$, $\Gamma_{j0}^{(y)}$ are known, since they depend on the known
functions $\psi_{10}$, ..., $\psi_{N-1,0}$, while $\Gamma^{(x,N)}$, $\Gamma^{(y,N)}$ have been introduced as in \equ{inv}.
Choosing $\tilde r_0<r_0$, $\delta_0<s_0$, the inversion is possible provided that (see Appendix~A)
\beq{33ter}
70\,\|\psi^{(N)}_y\|_{\tilde r_0,s_0}\ e^{2s_0}\delta_0^{-1}\ <\ 1\ ,
\eeq
being $\psi^{(N)}=\sum_{j=1}^N \psi_{j0}\varepsilon^j$.
For short, let us write the equations \equ{deltac3} as
$$
x=x(\tilde y,\tilde x,t)\ ,\qquad y=y(\tilde y,\tilde x,t)\ .
$$
Inserting \equ{deltac3} in  \equ{SE} and expanding in Taylor series, one has
\beqano
\dot x&=&\omega(\tilde y)+\omega_y(\tilde y)\psi_{N0,x}(\tilde y,\tilde x,t)\varepsilon^N+
F^{(0,N)}(\tilde y,\tilde x,t)\nonumber\\
&+&\varepsilon h^{>K}_{10,y}(y(\tilde y,\tilde x,t),
x(\tilde y,\tilde x,t),t)+\mu f^{>K}_{01}(y(\tilde y,\tilde x,t),x(\tilde y,\tilde x,t),t)
+O_{N+1}(\varepsilon,\mu)\noindent\\
\dot y&=&G^{(0,N)}(\tilde y,\tilde x,t)-\mu\left(\eta_{N-1,1}(\tilde y,\tilde x,t)\varepsilon^{N-1}+...+
\eta_{0,N}(\tilde y,\tilde x,t)\mu^{N-1}\right)\nonumber\\
&-&\varepsilon h^{>K}_{10,x}(y(\tilde y,\tilde x,t),x(\tilde y,\tilde x,t),t)+
\mu g^{>K}_{01}(y(\tilde y,\tilde x,t),x(\tilde y,\tilde x,t),t)
+O_{N+1}(\varepsilon,\mu)\ ,
\eeqano
where $F^{(0,N)}$, $G^{(0,N)}$ are known functions; $F^{(0,N)}$ contains terms of order
$\varepsilon$, $\varepsilon^2$, ..., $\varepsilon^N$,
$\mu$, $\mu\varepsilon$, ..., $\mu\varepsilon^{N-1}$, while $G^{(0,N)}$ contains all
powers $\varepsilon^j\mu^m$ with $1\leq j+m\leq N$. Next step is to compute
$\dot{\tilde x}$, $\dot{\tilde y}$ as a function of $\tilde x$, $\tilde y$. Taking into
account \equ{A} and that by the inductive hypothesis $\psi^{(1)}$, ..., $\psi^{(N-1)}$ make the equations in normal form
up to the order $N-1$, we obtain
\beqano
\dot{\tilde x}&=&\omega(\tilde y)+\sum_{j=1}^N \bar F^{(1,N)}_{j0}(\tilde y)\varepsilon^j+
\Big[\omega_y(\tilde y)\psi_{N0,x}(\tilde y,\tilde x,t)+\omega(\tilde y)\psi_{N0,yx}(\tilde y,\tilde x,t)
+\psi_{N0,yt}(\tilde y,\tilde x,t)\nonumber\\
&+&\tilde F^{(1,N,\leq K)}_{N0}(\tilde y,\tilde x,t)\Big]\varepsilon^N
+\mu \sum_{j=1}^{N-1} F^{(1,N)}_{j1}(\tilde y,\tilde x,t)\varepsilon^j
+\sum_{j=1}^{N} F_{j0}^{(1,N,>K)}(\tilde y,\tilde x,t)\varepsilon^j\nonumber\\
&+&\varepsilon h^{>K}_{10,y}(y(\tilde y,\tilde x,t),x(\tilde y,\tilde x,t),t)+
\mu f^{>K}_{01}(y(\tilde y,\tilde x,t),x(\tilde y,\tilde x,t),t)+O_{N+1}(\varepsilon,\mu)\ ,
\eeqano
where $F^{(1,N)}(\tilde y,\tilde x,t)$ is a known function that has been decomposed as
\beqano
F^{(1,N)}(\tilde y,\tilde x,t)&\equiv& \sum_{j=1}^N \bar F^{(1,N)}_{j0}(\tilde y)\varepsilon^j+
\sum_{j=1}^{N-1}\tilde F^{(1,N)}_{j0}(\tilde y,\tilde x,t)\varepsilon^j+
\tilde F^{(1,N)}_{N0}(\tilde y,\tilde x,t)\varepsilon^N\nonumber\\
&+&\mu \sum_{j=1}^{N-1} F^{(1,N)}_{j1}(\tilde y,\tilde x,t)\varepsilon^j\ .
\eeqano
In a similar way one obtains
\beqano
\dot{\tilde y}&=&\mu F^{(2,N)}(\tilde y,\tilde x,t)+\Big(\tilde F^{(3,N,\leq K)}_{N0}(\tilde y,\tilde x,t)-\omega(\tilde y)\psi_{N0,xx}(\tilde y,\tilde x,t)-
\psi_{N0,xt}(\tilde y,\tilde x,t)\Big)\varepsilon^N\nonumber\\
&-&\mu\left(\eta_{N-1,1}(\tilde y,\tilde x,t)\varepsilon^{N-1}+...+\eta_{0,N}(\tilde y,\tilde x,t)\mu^{N-1}\right)
+\sum_{j=1}^N F_{j0}^{(3,N,>K)}(\tilde y,\tilde x,t)\varepsilon^j\nonumber\\
&-&\varepsilon h^{>K}_{10,x}(y(\tilde y,\tilde x,t),x(\tilde y,\tilde x,t),t)+
\mu g^{>K}_{01}(y(\tilde y,\tilde x,t),x(\tilde y,\tilde x,t),t)
+O_{N+1}(\varepsilon,\mu)\ ,
\eeqano
for known functions $F^{(2,N)}$, $\tilde F^{(3,N)}$, the latter having zero average.
The conservative normal form is obtained by imposing that
$\psi_{N0}$ solves the following normal form equations
\beqa{star3}
\omega_y(\tilde y)\psi_{N0,x}(\tilde y,\tilde x,t)+\omega(\tilde y)
\psi_{N0,yx}(\tilde y,\tilde x,t)+\psi_{N0,yt}(\tilde y,\tilde x,t)+
\tilde F^{(1,N,\leq K)}_{N0}(\tilde y,\tilde x,t)&=&0\nonumber\\
\omega(\tilde y)\psi_{N0,xx}(\tilde y,\tilde x,t)+\psi_{N0,xt}(\tilde y,\tilde x,t)-
\tilde F^{(3,N,\leq K)}_{N0}(\tilde y,\tilde x,t)&=&0\ ,
\eeqa
where as before we have split the known function $F_{N0}^{(1,N)}$ into $F_{N0}^{(1,N)}=$ $\bar F_{N0}^{(1,N)}+\tilde F_{N0}^{(1,N)}$ as well as $\tilde F_{N0}^{(1,N)}$ into $\tilde F_{N0}^{(1,N)}=\tilde F_{N0}^{(1,N,\le K)}+ F_{N0}^{(1,N,>K)}$. Note that in this setting the $N$--th order contribution to the shifted frequency vector is given by $\Omega_{N0}(\tilde y)\equiv\bar F^{(1,N)}_{N0}(\tilde y)$. Due to the Hamiltonian character which occurs for $\mu=0$, there exists a function $M^{(N,\leq K)}=M^{(N,\leq K)}(\tilde y,\tilde x,t)$ with zero average, such that
$$
{{\partial M^{(N,\leq K)}(\tilde y,\tilde x,t)}\over {\partial y}}=\tilde
F^{(1,N,\leq K)}_{N0}(\tilde y,\tilde x,t)\ ,\qquad
{{\partial M^{(N,\leq K)}(\tilde y,\tilde x,t)}\over {\partial x}}=-\tilde
F^{(3,N,\leq K)}_{N0}(\tilde y,\tilde x,t)\ ,
$$
so that \equ{star3} are equivalent to solve the equation
\beq{psin0}
\omega(\tilde y)\psi_{N0,x}(\tilde y,\tilde x,t)+\psi_{N0,t}(\tilde y,\tilde x,t)
+M^{(N,\leq K)}(\tilde y,\tilde x,t)=0\ .
\eeq
The solution of \equ{psin0} provides the function $\psi_{N0}(\tilde y,\tilde x,t)$,
which produces the conservative normal form:
\beqa{CNF}
\dot{\tilde x}&=&\Omega_c^{(N)}(\tilde y;\varepsilon)+\mu \sum_{j=1}^{N-1} F_{j1}^{(1,N)}(\tilde y,\tilde x,t)\varepsilon^j
+\sum_{j=1}^N F_{j0}^{(1,N,>K)}(\tilde y,\tilde x,t)\varepsilon^j\nonumber\\
&+&\varepsilon h^{>K}_{10,y}(y(\tilde y,\tilde x,t),x(\tilde y,\tilde x,t),t)+
\mu f^{>K}_{01}(y(\tilde y,\tilde x,t),x(\tilde y,\tilde x,t),t)
+O_{N+1}(\varepsilon,\mu)\nonumber\\
\dot{\tilde y}&=&\mu F^{(2,N)}(\tilde y,\tilde x,t)-\mu\left(\eta_{N-1,1}(\tilde y,\tilde x,t)
\varepsilon^{N-1}+...+\eta_{0,N}(\tilde y,\tilde x,t)\mu^{N-1}\right)\nonumber\\
&-&\varepsilon h^{>K}_{10,x}(y(\tilde y,\tilde x,t),x(\tilde y,\tilde x,t),t)+
\mu g^{>K}_{01}(y(\tilde y,\tilde x,t),x(\tilde y,\tilde x,t),t)\nonumber\\
&+&\sum_{j=1}^N F_{j0}^{(3,N,>K)}(\tilde y,\tilde x,t)\varepsilon^j+O_{N+1}(\varepsilon,\mu)\ ,
\eeqa
where
\beq{omegaN}
\Omega_c^{(N)}(\tilde y;\varepsilon)\equiv\omega(\tilde y)+\sum_{j=1}^N\bar F_{j0}^{(1,N)}(\tilde y)
\varepsilon^j\ .
\eeq
Notice that from the definition \equ{omegaN} one obtains an estimate like
$\|\Omega_c^{(N)}-\omega\|\leq C_c\varepsilon$ for a suitable constant $C_c$.
The solution of \equ{psin0} involves small divisors of the form $\omega(\tilde y)\cdot k+m$,
for $k\in{\Z}^\ell$, $m\in{\Z}$ with $|k|+|m|\leq K$, as it was for the case $N=1$; non--resonance is
guaranteed provided that (see Appendix~A)
\beq{C6}
\varepsilon\leq{a\over {2K\|R^{(N)}\|_{\tilde r_0',s_0}\|\omega_y\|_{r_0}}}\ ,
\eeq
where $R^{(N)}$ is the function such that $\tilde y=y+\varepsilon R^{(N)}(y,x,t)$ and $\tilde r_0'<\tilde r_0$.

\vskip.1in

\noindent
\bf Step 4: Dissipative transformation for the order $N$. \rm

\noindent
As for the dissipative part, we consider the transformation \equ{B}, that we invert as
\beqa{star5}
\tilde x&=&X+\sum_{k=0}^N\sum_{j=0}^{N-k} a_{kj}(Y,X,t)\varepsilon^k\mu^j-
\sum_{k=0}^{N-1} \alpha_{k,N-k}(Y,X,t)\varepsilon^k\mu^{N-k}+O_{N+1}(\varepsilon,\mu)\nonumber\\
\tilde y&=&Y+\sum_{k=0}^N\sum_{j=0}^{N-k} b_{kj}(Y,X,t)\varepsilon^k\mu^j-
\sum_{k=0}^{N-1} \beta_{k,N-k}(Y,X,t)\varepsilon^k\mu^{N-k}+O_{N+1}(\varepsilon,\mu)\ ,
\eeqa
for suitable known functions $a_{kj}(Y,X,t)$, $b_{kj}(Y,X,t)$ with
$a_{k0}=b_{k0}=0$ for $k=0,...,N$, provided the following conditions are satisfied (see Appendix~A):
\beqa{C7}
70\,\|\alpha^{(N)}\|_{\tilde r_0,\tilde s_0}\ e^{2\tilde s_0}\tilde \delta_0^{-1}&<&1\nonumber\\
70\,\left(\|\beta^{(N)}\|_{\tilde r_0,\tilde s_0}+\|\beta_x^{(N)}\|_{\tilde r_0,\tilde s_0}\|\alpha^{(N)}\|_{\tilde r_0,\tilde s_0}\right)\ {1\over {\tilde r_0-R_0}}&<&1\ ,
\eeqa
where $\tilde\delta_0\equiv \tilde s_0/2$, $R_0<\tilde r_0$.
For short let us denote \equ{star5} as
$$
\tilde x=\tilde x(Y,X,t)\ ,\qquad \tilde y=\tilde y(Y,X,t)\ ,
$$
while we express the original variables in terms of $(Y,X,t)$ through
\beqano
x&=&x(\tilde y(Y,X,t),\tilde x(Y,X,t),t)\equiv X+\Phi^{(x,N)}(Y,X,t)\nonumber\\
y&=&y(\tilde y(Y,X,t),\tilde x(Y,X,t),t)\equiv Y+\Phi^{(y,N)}(Y,X,t)\ ,
\eeqano
for suitable functions $\Phi^{(x,N)}$, $\Phi^{(y,N)}$, which are $O_1(\varepsilon,\mu)$.
We need to determine the unknown
functions $\alpha_{0,N}$, ..., $\alpha_{N-1,1}$, $\beta_{0,N}$, ..., $\beta_{N-1,1}$,
$\eta_{N-1,1}$, ..., $\eta_{0,N}$ as follows. Starting from \equ{CNF}, we compute
$\dot{\tilde x}$, $\dot{\tilde y}$ in terms of $X$, $Y$ and we express $\dot X$, $\dot Y$ in terms of $X$, $Y$,
using \equ{B}, \equ{CNF}, \equ{star5}; by the inductive hypothesis $\alpha_{kj}$, $\beta_{kj}$, $\eta_j$
with $0\leq k+j\leq N-1$ are determined so that the equations of motion are in normal form
up to the order $\varepsilon^k\mu^j$ with $0\leq k+j\leq N-1$. This leads to the following equations:
\beqano
\dot X&=&\Omega_c^{(N)}(Y)-\omega_y(Y)\left(\sum_{j=0}^{N-1}\beta_{j,N-j}(Y,X,t)\varepsilon^j\mu^{N-j}\right)
+\omega(Y)\sum_{j=0}^{N-1}\alpha_{j,N-j,x}(Y,X,t)\varepsilon^j\mu^{N-j}\nonumber\\
&+&\sum_{j=0}^{N-1}\alpha_{j,N-j,t}(Y,X,t)\varepsilon^j\mu^{N-j}
+\mu F^{(4,N)}(Y,X,t)\nonumber\\
&+&\sum_{j=1}^N F_{j0}^{(1,N,>K)}(\tilde y(Y,X,t),\tilde x(Y,X,t),t)\varepsilon^j\nonumber\\
&+&\varepsilon h^{>K}_{10,y}(y(Y,X,t),x(Y,X,t),t)+\mu f^{>K}_{01}(y(Y,X,t),x(Y,X,t),t)
+O_{N+1}(\varepsilon,\mu)\nonumber\\
\dot Y&=&-\mu\Big(\eta_{N-1,1}(Y,X,t)\varepsilon^{N-1}+ ...+\eta_{0,N}(Y,X,t)\mu^{N-1}\Big)
+\omega(Y)\sum_{j=0}^{N-1}\beta_{j,N-j,x}(Y,X,t)\varepsilon^j\mu^{N-j}\nonumber\\
&+&\sum_{j=0}^{N-1}\beta_{j,N-j,t}(Y,X,t)\varepsilon^j\mu^{N-j}
+\mu F^{(5,N)}(Y,X,t)\nonumber\\
&+&\sum_{j=1}^N F_{j0}^{(3,N,>K)}(\tilde y(Y,X,t),\tilde x(Y,X,t),t)\varepsilon^j\nonumber\\
&-&\varepsilon h^{>K}_{10,x}(y(Y,X,t),x(Y,X,t),t)+\mu g^{>K}_{01}(y(Y,X,t),x(Y,X,t),t)
+O_{N+1}(\varepsilon,\mu)\ ,
\eeqano
where $F^{(4,N)}(Y,X,t)$, $F^{(5,N)}(Y,X,t)$ are $O(\mu^{N-1},\varepsilon\mu^{N-2},...,\varepsilon^{N-2}\mu,
\varepsilon^{N-1})$. Let us decompose
$F^{(4,N)}$ as $F^{(4,N)}(Y,X,t)=\bar F^{(4,N)}(Y)+\tilde F^{(4,N)}(Y,X,t)$. The normal form
at the order $N$ is achieved imposing that $\alpha_{kj}$, $\beta_{kj}$, $\eta_j$ satisfy
the following normal form equations
\beqa{star7}
&&-\omega_y(Y)\beta_{j,N-j}(Y,X,t)+\omega(Y)\alpha_{j,N-j,x}(Y,X,t)+
\alpha_{j,N-j,t}(Y,X,t)\nonumber\\
&&\qquad\qquad\qquad\qquad\qquad\qquad\qquad\qquad+\tilde F^{(4,N,\leq K)}_{j,N-j}(Y,X,t)=0\nonumber\\
&&\omega(Y)\beta_{j,N-j,x}(Y,X,t)+
\beta_{j,N-j,t}(Y,X,t)+\tilde F^{(5,N,\leq K)}_{j,N-j}(Y,X,t)\nonumber\\
&&\qquad\qquad\qquad\qquad\qquad\qquad\qquad\qquad+\bar F_{j,N-j}^{(5,N)}(Y)-\eta_{j,N-j}(Y,X,t)=0
\eeqa
for $0\leq j\leq N-1$, where we have used the expansion
$$
\mu F^{(m,N)}(Y,X,t)\equiv
\sum_{k=0}^{N-1} \bar F^{(m,N)}_{k,N-k}(Y) \varepsilon^k\mu^{N-k}
+\sum_{k=0}^{N-1} \tilde F^{(m,N)}_{k,N-k}(Y,X,t)
\varepsilon^k\mu^{N-k}\ ,\quad m=4,\ 5
$$
and we split $\tilde F^{(m,N)}_{k,N-k}=\tilde F^{(m,N,\leq K)}_{k,N-k}+F^{(m,N,>K)}_{k,N-k}$.
The non--resonance condition for $\omega(Y)$ reads as (see Appendix~A):
\beq{C8}
K\|\ \beta^{(N)}\|_{\tilde r_0,\tilde s_0} \|\omega_y\|_{r_0}<{a\over 4}\ ,
\eeq
where, for short, we have written $Y\equiv \tilde y+\beta^{(N)}(\tilde y,\tilde x,t;\varepsilon,\mu)$.
From the second of \equ{star7}, it is
$$
\eta_{j,N-j}(Y,X,t)\equiv \eta_{j,N-j}(Y)=\bar F^{(5,N)}_{j,N-j}(Y)\ ,
$$
so that the second of \equ{star7} can be solved to determine $\beta_{j,N-j}$, while
from the first of \equ{star7} we obtain $\alpha_{j,N-j}$. Identifying the final normalized frequency
$\Omega_d$ with $\Omega_d^{(N)}(Y;\varepsilon,\mu)\equiv \omega(Y)+\sum_{m=0}^N\sum_{j=0}^m\Omega_{j,m-j}(Y)\varepsilon^j\mu^{m-j}\ ,$
where $\Omega_{00}=0$, $\Omega_{10}=\omega$, $\Omega_{j0}\equiv \bar F_{j0}^{(1,N)}(Y)$ and $\Omega_{j,m-j}\equiv \bar F^{(4,N)}_{j,m-j}(Y)$, the resulting normal form is given by
\beqa{NFf}
\dot X &=& \Omega_d^{(N)}(Y;\varepsilon,\mu)+F_{N+1}(Y,X,t)
+\sum_{j=0}^{N-1} F^{(4,N,>K)}_{j,N-j}(Y,X,t)\varepsilon^j\mu^{N-j}\nonumber\\
&+&\varepsilon h^{>K}_{10,y}(y(Y,X,t),x(Y,X,t),t)+\mu f^{>K}_{01}(y(Y,X,t),x(Y,X,t),t)\nonumber\\
&+&\sum_{j=1}^N F^{(1,N,>K)}_{j0}(\tilde y(Y,X,t),\tilde x(Y,X,t),t)\varepsilon^j\nonumber\\
\dot Y  &=& G_{N+1}(Y,X,t)+\sum_{j=0}^{N-1} F^{(5,N,>K)}_{j,N-j}(Y,X,t)\varepsilon^j\mu^{N-j}
\nonumber\\&-&\varepsilon h^{>K}_{10,x}(y(Y,X,t),x(Y,X,t),t)+\mu g^{>K}_{01}(y(Y,X,t),x(Y,X,t),t)\nonumber\\
&+&\sum_{j=1}^N F^{(3,N,>K)}_{j0}(\tilde y(Y,X,t),\tilde x(Y,X,t),t)\varepsilon^j\ ,
\eeqa
where $F_{N+1}$, $G_{N+1}$ are $O_{N+1}(\varepsilon,\mu)$. For short, we define
\beqano
F_{01}^{>K}(Y,X,t)&\equiv& \sum_{j=0}^{N-1} F^{(4,N,>K)}_{j,N-j}(Y,X,t)\varepsilon^j\mu^{N-j}
+\sum_{j=1}^N F_{j0}^{(1,N,>K)}(\tilde y(Y,X,t),\tilde x(Y,X,t),t)\varepsilon^j\nonumber\\
G_{01}^{>K}(Y,X,t)&\equiv& \sum_{j=0}^{N-1} F^{(5,N,>K)}_{j,N-j}(Y,X,t)\varepsilon^j\mu^{N-j}
+\sum_{j=1}^N F_{j0}^{(3,N,>K)}(\tilde y(Y,X,t),\tilde x(Y,X,t),t)\varepsilon^j\ ,
\eeqano
which makes \equ{NFf} of the form \equ{NF}. The estimate \equ{Pi} follows from the fact that \equ{A} is close to the identity up to terms of order $\varepsilon$, while \equ{B} is close
to the identity up to terms of order $\mu$. Notice that the smallness requirements on $\varepsilon$, $\mu$, say $\varepsilon\leq \varepsilon_0$, $|\mu|\leq\mu_0$, are needed to ensure that the non--resonance condition (see \equ{C3}, \equ{C5}, \equ{C6}, \equ{C8}) is satisfied and that the transformations \equ{A}, \equ{B} can be inverted (see \equ{C1}, \equ{C2}, \equ{C4}, \equ{33ter}, \equ{C7}).

\vskip.1in

\noindent
The original variables can be expressed in terms of the intermediate variables by means of \equ{deltac3}, provided
\equ{33ter} holds with
$$
\|\Gamma^{(x,N)}\|_{\tilde r_0,\tilde s_0}\leq \|\psi_y^{(N)}\|_{\tilde r_0,s_0}\ ,
$$
being $\tilde s_0\equiv s_0-\delta_0$. Moreover, we have
$$
\|\Gamma^{(y,N)}\|_{\tilde r_0,\tilde s_0}\leq \|\psi_x^{(N)}\|_{\tilde r_0,s_0}+
\|\psi_{xx}^{(N)}\|_{\tilde r_0,s_0}\|\Gamma^{(x,N)}\|_{\tilde r_0,\tilde s_0}\ .
$$
For the same reason, equations \equ{star5} are invertible for $\varepsilon$, $|\mu|$ sufficiently small, since
the Jacobian of the transformation is close to the identity; we can write the inverse as
\beqano
\tilde x&=& X+\Delta^{(x,N)}(Y,X,t)\nonumber\\
\tilde y&=& Y+\Delta^{(y,N)}(Y,X,t)
\eeqano
for suitable, bounded functions $\Delta^{(x,N)}$, $\Delta^{(y,N)}$ of order $O_1(\varepsilon,\mu)$. We finally obtain
\beqano
x&=& X+\Delta^{(x,N)}(Y,X,t)+\Gamma^{(x,N)}(Y+\Delta^{(y,N)}(Y,X,t),X+\Delta^{(x,N)}(Y,X,t),t)\nonumber\\
&\equiv& X+\Phi^{(x,N)}(Y,X,t)\nonumber\\
y&=& Y+\Delta^{(y,N)}(Y,X,t)+\Gamma^{(y,N)}(Y+\Delta^{(y,N)}(Y,X,t),X+\Delta^{(x,N)}(Y,X,t),t)\nonumber\\
&\equiv& Y+\Phi^{(y,N)}(Y,X,t)\ .
\eeqano
Recalling Lemma B.1 of Appendix~B, for $\tau_0>0$ and $1\leq j\leq \ell$, we have
\beqano
\|h^{>K}_{10,x}\|_{R_0,S_0}&\leq& C_a \|h_{10,x}\|_{R_0,S_0+\tau_0}e^{-(K+1)\tau_0}\nonumber\\
\|h^{>K}_{10,xy_j}\|_{R_0,S_0}&\leq& C_a \|h_{10,xy_j}\|_{R_0,S_0+\tau_0}e^{-(K+1)\tau_0}\nonumber\\
\|h^{>K}_{10,xx_j}\|_{R_0,S_0}&\leq& C_a \|h_{10,xx_j}\|_{R_0,S_0+\tau_0}e^{-(K+1)\tau_0}\ ,
\eeqano
where $C_a$ is a positive real constant.
Setting
\beqano
C_{(\Phi,y)}&\equiv& \sup_{1\leq j\leq \ell}\|\Phi_j^{(y,N)}\|_{R_0,S_0}\ ,\qquad
C_{(\Phi,x)}\equiv \sup_{1\leq j\leq \ell}\|\Phi_j^{(x,N)}\|_{R_0,S_0}\ ,\nonumber\\
C_{(h,y)}&\equiv& \sup_{1\leq j\leq \ell}\|h_{10,xy_j}\|_{R_0,S_0}\ ,\qquad
C_{(h,x)}\equiv \sup_{1\leq j\leq \ell}\|h_{10,xx_j}\|_{R_0,S_0}\ ,
\eeqano
one obtains
\beqano
&&\|h^{>K}_{10,x}(Y+\Phi^{(y,N)}(Y,X,t),X+\Phi^{(x,N)}(Y,X,t),t)\|_{R_0,S_0}\leq
\|h^{>K}_{10,x}\|_{R_0,S_0}\nonumber\\
&+&\ell\Big(\sup_{1\leq j\leq \ell}\|h_{10,xy_j}\|_{R_0,S_0}
\sup_{1\leq j\leq \ell}\|\Phi_j^{(y,N)}\|_{R_0,S_0}
+\sup_{1\leq j\leq \ell}\|h_{10,xx_j}\|_{R_0,S_0}
\sup_{1\leq j\leq \ell}\|\Phi_j^{(x,N)}\|_{R_0,S_0}\Big)\nonumber\\
&\equiv& C_he^{-K\tau_0}\ ,
\eeqano
having defined $C_h\equiv C_a e^{-\tau_0}
\Big[\|h_{10,x}\|_{R_0,S_0+\tau_0}+\ell\Big(C_{(\Phi,y)}C_{(h,y)}+C_{(\Phi,x)}C_{(h,x)}\Big)\Big]$.
Similarly we obtain
\beqno
\|g^{>K}_{01}(Y+\Phi^{(y,N)}(Y,X,t),X+\Phi^{(x,N)}(Y,X,t),t)\|_{R_0,S_0}\nonumber\\
\leq C_ge^{-K\tau_0}\ ,
\eeqno
where
\beqano
C_{(g,y)}&\equiv& \sup_{1\leq j\leq \ell}\|g_{01,y_j}\|_{R_0,S_0}\ ,\qquad
C_{(g,x)}\equiv \sup_{1\leq j\leq \ell}\|g_{01,x_j}\|_{R_0,S_0}\ ,\nonumber\\
C_g&\equiv& C_a e^{-\tau_0}
\Big[\|g_{01}\|_{R_0,S_0+\tau_0}+\ell\Big(C_{(\Phi,y)}C_{(g,y)}+C_{(\Phi,x)}C_{(g,x)}\Big)\Big]\ .
\eeqano
Analogously we find
\beqno
\|G^{>K}_{01}(Y+\Phi^{(y,N)}(Y,X,t),X+\Phi^{(x,N)}(Y,X,t),t)\|_{R_0,S_0}\leq \lambda \tilde C_G e^{-K\tau_0}\ ,
\eeqno
where
\beqano
C_{(G,y)}&\equiv& \sup_{1\leq j\leq \ell}\|\mu_0^{-1} G_{01,y_j}^{>K}\|_{R_0,S_0}\ ,\qquad
C_{(G,x)}\equiv \sup_{1\leq j\leq \ell}\|\mu_0^{-1} G_{01,x_j}^{>K}\|_{R_0,S_0}\ ,\nonumber\\
\tilde C_G&\equiv& C_a e^{-\tau_0}
\Big[\|\mu_0^{-1} G_{01}^{>K}\|_{R_0,S_0+\tau_0}+\ell\Big(C_{(\Phi,y)}C_{(G,y)}+C_{(\Phi,x)}C_{(G,x)}\Big)\Big]\ .
\eeqano
Let us bound $G_{N+1}$ in \equ{NFf} as
$$
\|G_{N+1}\|_{R_0,S_0}\leq C_G\lambda^{N+1}\ ,
$$
for a suitable constant $C_G$. From the second of \equ{NFf} we obtain:
\beqano
\|G_{N+1}\|_{R_0,S_0}&+&\varepsilon \|h^{>K}_{10,x}\|_{R_0,S_0}+|\mu| \|g^{>K}_{01}\|_{R_0,S_0}
+\|G^{>K}_{01}\|_{R_0,S_0}\leq C_G\lambda^{N+1}\nonumber\\
&+&\varepsilon C_he^{-K\tau_0}+|\mu| C_ge^{-K\tau_0}+\lambda \tilde C_G e^{-K\tau_0}\ .
\eeqano
Choosing $N$ as
\beq{tau0}
N\equiv [{{K\tau_0}\over |\log \lambda|}]\ ,
\eeq
we have (see \equ{FG})
$$
\|G_{N+1}\|_{R_0,S_0}+\varepsilon \|h^{>K}_{10,x}\|_{R_0,S_0}+|\mu| \|g^{>K}_{01}\|_{R_0,S_0}
+\|G^{>K}_{01}\|_{R_0,S_0}\leq C_Y\lambda^{N+1}\ ,
$$
with $C_Y\equiv C_G+C_h+C_g+\tilde C_G$. This concludes the proof of the Lemma. $\Box$

\vskip.2in

\noindent
\bf Proof of the theorem. \rm The distance between the solution at time $t>0$, say $y(t)$, and the initial
condition $y(0)$ can be bounded as
$$
\|y(t)-y(0)\|\leq \|y(t)-Y(t)\|+\|Y(t)-Y(0)\|+\|Y(0)-y(0)\|\ .
$$
By the estimate \equ{Pi} of the Normal Form Lemma one has
$$
\|y(t)-Y(t)\|\leq C_p \lambda\ ,\qquad \|y(0)-Y(0)\|\leq C_p \lambda\ .
$$
By the second of \equ{NF} and by \equ{FG}, one has
\beqano
\|Y(t)-Y(0)\|&\leq& \int_0^t
\left(\|G_{N+1}\|_{R_0,S_0}+\varepsilon \|h^{>K}_{10,x}\|_{R_0,S_0}
+|\mu| \|g^{>K}_{01}\|_{R_0,S_0}+\|G^{>K}_{01}\|_{R_0,S_0}\right)ds\nonumber\\
&\leq& t C_Y\lambda^{N+1}\ .
\eeqano
Let $r_2>0$ be such that
$$
t\  C_Y\lambda^{N}\leq r_2\ ,
$$
which is satisfied as far as
$$
t\leq {r_2\over C_Y}\ \lambda^{-N}\leq {r_2\over {C_Y}}\ e^{K\tau_0}\ ,
$$
where \equ{tau0} has been used. Finally, setting $\rho_0\equiv 2C_p+r_2$, we obtain
$$
\|y(t)-y(0)\|\leq \rho_0\lambda\qquad {\rm for\ }\ t\leq C_t e^{K\tau_0}\ ,
$$
having defined $C_t\equiv r_2/C_Y$. $\Box$

\section{An application of the normal form and of the stability estimates}\label{sec:example}
In order to test the accuracy of our results, we implement the Normal Form Lemma and we derive the stability estimates on a specific example. To this end, let $\ell=1$ and let us consider the differential system:
\beqa{e19}
\dot x&=&y-\mu \left(\sin (x-t)+\sin (x)\right) , \nonumber \\
\dot y&=&-\varepsilon \left(\sin (x-t)+\sin (x)\right) -\mu  (y-\eta)\ .
\eeqa
We remark that for $\mu=0$ the system \equ{e19} is associated to the Hamiltonian function in the
extended phase space
$$
H(y,T,x,t)=\frac{y^2}{2}-\varepsilon \left(\cos (x-t)+\cos (x)\right)+T \ ,
$$
where the unperturbed frequency of the motion is given by $\omega(y)=y$ and being $T$ conjugated to time.
We provide details for the computation of the second order normal form associated to \equ{e19} (see Section~\ref{sec:nf}).
A comparison with a numerical integration is performed in Section~\ref{sec:comparison}.
Stability estimates according to the Theorem of Section~\ref{sec:stability} are computed in Section~\ref{sec:estimates}. A slightly different example with oscillating energy is analyzed in Section~\ref{sec:oscillating}.

\subsection{Normal form}\label{sec:nf}
The second order normal form can be computed as follows \footnote{Notice that here we first implement the conservative transformation to the second order and then we determine the dissipative transformation, which provides the second order normal form.}. At first order we identify the non--zero average contributions by
\(\bar h_{10,y}(\tilde{y})=0\) (see \equ{omegag}). The conservative normal form equations become (see \equ{otto}):
\beqano
\psi _{10,x}(\tilde y,\tilde x,t)+ \psi _{10,yt}(\tilde y,\tilde x,t)+\tilde{y} \ \psi _{10,yx}(\tilde y,\tilde x,t) &=&0 \nonumber \\
\tilde{y} \ \psi _{10,xx}(\tilde y,\tilde x,t)+\psi _{10,xt}(\tilde y,\tilde x,t)
+\sin (\tilde{x}-t)+\sin (\tilde{x}) &=& 0 \ ,
\eeqano
from which we get
\beqno
\psi _{10}(\tilde{y},\tilde{x},t)=\frac{\sin (\tilde{x}-t)}{\tilde{y}-1}+\frac{\sin (\tilde{x})}{\tilde{y}} \ .
\eeqno
The second order conservative normal form is obtained by computing the generating function $\psi_{20}$
as the solution of the equations:
\beqano
\psi _{20,x}(\tilde y,\tilde x,t)+\psi _{20,yt}(\tilde y,\tilde x,t)+\tilde{y} \ \psi _{20,yx}(\tilde y,\tilde x,t)+\frac{\cos (2 \tilde{x}-2t)}{2 (1-\tilde{y})^3}&&\nonumber\\
+\frac{(1-2 \tilde{y} ) \cos (2 \tilde{x}-t)}{2 (\tilde{y}-1)^2 \tilde{y}^2}+\frac{\cos (t) (1-2 \tilde{y})}{2
(\tilde{y}-1)^2 \tilde{y}^2}-\frac{\cos (2 \tilde{x})}{2 \tilde{y}^3} &=& 0, \nonumber\\
\tilde{y} \ \psi _{20,xx}(\tilde y,\tilde x,t)+\psi _{20,xt}(\tilde y,\tilde x,t)+\frac{\sin (2 \tilde{x}-t)}{(\tilde{y}-1) \tilde{y}}+\frac{\sin (2 \tilde{x}-2t)}{2 (\tilde{y}-1)^2}+\frac{\sin (2 \tilde{x})}{2 \tilde{y}^2} &=& 0 \ ,
\eeqano
\noindent
which provides the function \(\psi _{20}\) as
\beqno
\psi_{20}(\tilde{y},\tilde{x},t)=-\frac{\sin (2 \tilde{x}-t)}{2 (\tilde{y}-1) \tilde{y} (2 \tilde{y}-1)}-\frac{\sin (2 \tilde{x}-2t)}{8 (\tilde{y}-1)^3}-\frac{\sin (t)}{2 (\tilde{y}-1) \tilde{y}}-\frac{\sin (2 \tilde{x})}{8 \tilde{y}^3} \ ,
\eeqno
while the second order term of the frequency shift is given by
\(\Omega_{20}(\tilde{y})=(-2 \tilde{y}^3+3 \tilde{y}^2 -3 \tilde{y} +1)/(2 (\tilde{y}-1)^3 \tilde{y}^3)\), so that
$\Omega_c^{(2)}(\tilde y)=\omega(\tilde y)+\Omega_{20}(\tilde y)\varepsilon^2$ (see \equ{omegaN}). At this stage, we succeeded in normalizing the symplectic contributions and in getting the conservative normal form to the second order in the
intermediate variables as
\beqano
\dot{\tilde{x}}&=&\tilde{y}+\frac{\left(1-3 \tilde{y}+3 \tilde{y}^2-2 \tilde{y}^3\right) }{2 (\tilde{y}-1)^3 \tilde{y}^3}\varepsilon ^2-\mu(\sin(\tilde x-t)+\sin(\tilde x))+O_3(\varepsilon) \nonumber \\
\dot{\tilde{y}}&=&-\mu(\tilde y-\eta)+O_3(\varepsilon)\ .
\eeqano
The first order of the dissipative normal form, expressed in terms of the new variables, provides the equations
\beqano
Y \alpha_{01,x}(Y,X,t)+\alpha _{01,t}(Y,X,t)-\beta _{01}(Y,X,t)-\sin (X-t)-\sin (X) &=& 0 \nonumber \\
Y \beta_{01,x}(Y,X,t)+\beta _{01,t}(Y,X,t)+Y-\eta_{01}(Y,X,t) &=& 0
\eeqano
(see \equ{XX}). From the second equation we get
$$
\eta _{01}(Y,X,t)=Y\ ,\qquad \beta _{01}(Y,X,t)=0 \ ,
$$
while from the first equation we obtain
\beqano
\alpha _{01}(Y,X,t)&=&\frac{\cos (X-t)}{1-Y}-\frac{\cos (X)}{Y} \ .
\eeqano
In a similar way, the second order dissipative normal form provides the equations
\beqa{alfabeta}
&&Y \alpha_{11,x}(Y,X,t)+\alpha_{11,t}(Y,X,t)-\beta_{11}(Y,X,t)+\frac{(2 Y-1) \sin (t)}{(Y-1)^2 Y^2} = 0, \nonumber \\
&&Y \beta_{11,x}(Y,X,t)+\beta_{11,t}(Y,X,t)+\eta _{11}(Y,X,t)+\frac{1-2 Y}{2 (Y-1) Y}-\frac{(2 Y-1) \cos (2 X-t)}{2 (Y-1) Y}- \nonumber \\
&&\qquad\qquad\qquad \frac{\cos (2 X-2 t)}{2 (Y-1)}+\frac{(1-2 Y) \cos (t)}{2 (Y-1) Y}-\frac{\cos (2 X)}{2 Y} = 0 \nonumber\\
&&Y \alpha_{02,x}(Y,X,t)+\alpha_{02,t}(Y,X,t)-\beta _{02}(Y,X,t)+\frac{(2 Y-1) \cos (2 X-t)}{2 (Y-1) Y}+\frac{\cos (2 X-2 t)}{2 (Y-1)}+ \nonumber \\
&&\qquad\qquad\qquad \frac{(1-2 Y) \cos (t)}{2 (Y-1) Y}+\frac{\cos (2 X)}{2 Y} = 0, \nonumber \\
&&Y \beta_{02,x}(Y,X,t)+\beta_{02,t}(Y,X,t)+\eta_{02}(Y,X,t) = 0\ .
\eeqa
From the second and fourth of \equ{alfabeta} we conclude that
\beqano
\eta_{02}(Y)&=&0 \nonumber \\
\eta_{11}(Y)&=&\frac{2 Y-1}{2 (Y-1) Y} \nonumber \\
\beta_{11}(Y,X,t)&=&\frac{\sin (2 X-t)}{2 (1-Y) Y}-
\frac{\sin (2 X-2 t)}{4 (Y-1)^2}-\frac{(1-2 Y) \sin (t)}{2 (Y-1) Y}-\frac{\sin (2 X)}{4 Y^2} \nonumber \\
\beta_{02}(Y,X,t)&=&0\ ,
\eeqano
while from the first and third of \equ{alfabeta} we obtain
\beqano
\alpha_{11}(Y,X,t)&=&\frac{\cos (2 X-t)}{2 Y \left(2 Y^2-3 Y+1\right)}+\frac{\cos (2 X-2 t)}{8 (Y-1)^3}+ \nonumber \\
&&\frac{\left(-2 Y^3+3 Y^2+3 Y-2\right) \cos(t)}{2 (Y-1)^2 Y^2}+\frac{\cos (2 X)}{8 Y^3} \nonumber \\
\alpha_{02}(Y,X,t)&=&\frac{\sin (2 X-t)}{2 (Y-1) Y}+\frac{\sin (2 X-2 t)}{4 (Y-1)^2}+\frac{(2 Y-1) \sin (t)}{2 (Y-1) Y}+\frac{\sin (2 X)}{4 Y^2} \ .
\eeqano
In conclusion, the second order normal form associated to \equ{e19} is given by
\beqano
\dot X&=&Y+\frac{\left(1-3Y+3Y^2-2Y^3\right) }{2 (Y-1)^3 Y^3}\varepsilon ^2+\frac{(1-2 Y) }{2 (Y-1) Y}\mu ^2+O_3(\varepsilon,\mu)  \nonumber \\
\dot Y&=&0+O_3(\varepsilon,\mu) \ ,
\eeqano
where $\eta=\eta(Y)$ takes the following expression:
\beqno
\eta(Y) =Y-\frac{1-2 Y}{2 (Y-1) Y}\varepsilon +O_3(\varepsilon,\mu)\ .
\eeqno
This concludes the computation for the second order normal form. In a similar way one can continue
to higher orders and in fact we computed up to the fifth order were we stopped, since i) we already reached exponential estimates for non trivial parameter values (see Section~\ref{sec:estimates}) and ii) the structure of the terms to be determined becomes too complex to allow for higher order computations using just a general purpose algebraic manipulator ({\sl Mathematica 7}); nevertheless, we believe that higher orders can be obtained by implementing a specific algebraic manipulator in C or Fortran languages.
In order to understand the degree of complexity of the computation,
let us denote by $\Xi_c^{(N)}$, $\Xi_d^{(N)}$, $\Xi_c^{(N)}\circ\Xi_d^{(N)}$ respectively, the conservative normal form at order $N$, the dissipative transformation and the overall normal form. At any order $N$ the algebraic manipulator has to deal with Poisson series (\cite{JHEN89}) of the form
$$
\sum _{(j,k)\in U\subset {\Z}^{\ell+1}} a_{jk}\varepsilon^{b_{jk}} \mu^{c_{jk}}\frac{P_{jk}(y)}{Q_{jk}(y)}e^{-i (j\cdot x+k t)} \ ,
$$
where $U$ is a sublattice of ${\Z}^{\ell+1}$,
$a_{jk}$ are complex coefficients, $b_{jk}, c_{jk}\in{\Z}_+$ with $b_{jk}+c_{jk}=N$ and $P_{jk}$, $Q_{jk}$ are polynomials in the actions. Let us denote by $deg(P_{jk})$ the degree of the polynomial $P_{jk}$, which contains only positive (or zero) powers in the action (and similarly for $Q_{jk}$).
For each order $N$ between 1 and 5 the numbers of Fourier terms as well as the degree
of the polynomials $P_{jk}$ and $Q_{jk}$ are provided in Table~\ref{effort}.
We remark that the main limitation in the present implementation of the normal form algorithm turned out to be the capability of dealing with the algebraic manipulation of fractions of polynomials of higher degree ($29$ at order $5$)\footnote{Note, that the number of terms of a Poisson series expansion with rational coefficients in the actions strongly depends on the 'algebraic normal form' of them. Putting terms of the same Fourier mode under the same denominator will increase the exponent order in the denominators and reduce the number of Fourier terms. On the contrary, writing the sums apart will reduce the exponent order in the denominator, but increase the number of terms. The numbers given in Table \ref{effort} strongly depend on the choice of the form of the rational coefficients.}.

\vskip.2in

\begin{table}
\caption{Number of Fourier terms and the degree of the polynomials of the conservative and dissipative transformations as a function of the order.}

\vskip.1in

\centering

\begin{tabular}{|c|c|c|c|c|c|}
  \hline
  & $N=1$ & $N=2$ & $N=3$ & $N=4$ & $N=5$ \\
  \hline
  $\#$ Fourier terms for $\Xi_c^{(N)}$ & 9 & 26 & 130 & 524 & 1888 \\
  $\#$ Fourier terms for $\Xi_d^{(N)}$ & 6 & 25 & 201 & 846 & 6829 \\
  $\#$ Fourier terms for $\Xi_c^{(N)}\circ\Xi_d^{(N)}$ & 11 & 53 & 461 & 2875 & 5004 \\
  $deg(P_{jk})$ for $\Xi_c^{(N)}$ & 0 & 4 & 6 & 8 & 10 \\
  $deg(P_{jk})$ for $\Xi_d^{(N)}$ & 1 & 3 & 5 & 7 & 9 \\
  $deg(P_{jk})$ for $\Xi_c^{(N)}\circ\Xi_d^{(N)}$ & 2 & 5 & 6 & 18 & 29 \\
  $deg(Q_{jk})$ for $\Xi_c^{(N)}$ & 2 & 1 & 1 & 1 & 2 \\
  $deg(Q_{jk})$ for $\Xi_d^{(N)}$ & 1 & 1 & 1 & 1 & 1 \\
  $deg(Q_{jk})$ for $\Xi_c^{(N)}\circ\Xi_d^{(N)}$ & 2 & 2 & 2 & 11 & 19 \\
  \hline
  \label{effort}
 \end{tabular}
\end{table}

\subsection{Comparison with a numerical integration}\label{sec:comparison}
To compare the results with a direct numerical integration we need to find the expression for $\eta(Y)$ in terms of the original variables \((y,x,t)\). To this end, we determine $\eta(y)$ from the condition that still $\dot Y=0$
in terms of the original variables up to the normalization order $N$, thus obtaining
\beqano
\eta(y)=y-\frac{(1-2y) }{2(y-1)y}\varepsilon +O_3(\varepsilon,\mu) \ .
\eeqano

\vskip.1in

\noindent
\bf Remark. \rm By induction on the normalization order, one can easily prove that the normal form equation
in new and old variables keeps the same functional form. For this reason also $\eta$ maintains the same
form in old and new variables.

\vskip.1in

\begin{figure}
\begin{center}
\includegraphics[width=7cm,keepaspectratio]{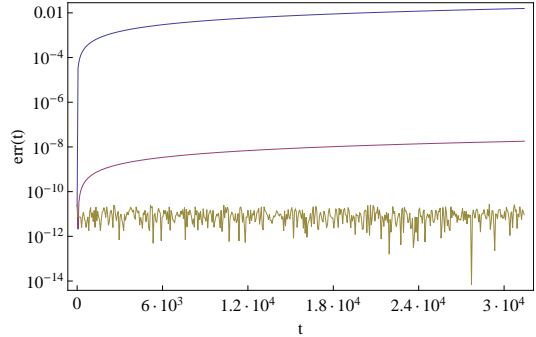}
\caption{Comparison between the analytical and numerical solutions obtained using normal forms of 1st (upper curve), 3rd (central curve) and 5th (bottom curve) orders. The integration time is \(T=10^4\pi \), the parameters are \(\varepsilon=10^{-3}\), \(\mu=10^{-3}\); the initial condition is set to \(Y_0 =\frac{1}{2}\left(\sqrt{5}+1\right)\), $X_0=0$.}
\label{fig-1}
\end{center}
\end{figure}

\noindent
For given initial conditions $(X_0,Y_0)$ we integrate the normal form equations (up to the normalization order) as
\beqano
X(t)&=&\Omega_d^{(N)} \left(Y_0; \varepsilon, \mu\right)t+X_0 \nonumber \\
Y(t)&=&Y_0 \ ,
\eeqano
where $\Omega_d^{(N)}=\Omega_d^{(N)}(Y;\varepsilon,\mu)$ denotes the normalized frequency to the $N$--th order. Then we back--transform to old variables and we define the relative error between the analytical and the numerical solution as
\beq{errt}
err(t)\equiv \frac{1}{2}\frac{\left(\left(x_{\rm{num}}(t)-x_{\rm{ana}}(t)\right){}^2+\left(y_{\rm{num}}(t)-y_{\rm{ana}}(t)
\right){}^2\right){}^{1/2}}{\left(x_{\rm{num}}(t){}^2+x_{\rm{ana}}(t){}^2+y_{\rm{num}}(t){}^2+y_{\rm{ana}}(t){}^2\right)
{}^{1/2}} \ ,
\eeq
where $\left(x_{\rm{num}},y_{\rm{num}}\right)$ is the state vector at time $t$ obtained from
a numerical integration of the original equations of motion and $\left(x_{\rm{ana}},y_{\rm{ana}}\right)$ represents the state vector at time \(t\) obtained from the normal form solution back--transformed in the original variables.
The evolution in time of the error for the parameter values \(\varepsilon=10^{-3}\),
$\mu =10^{-3}$ and the initial conditions \(Y_0=\frac{1}{2}(\sqrt{5}+1)\), $X_0=0$,
is shown in Figure~\ref{fig-1}. We plot the value of $err(t)$ versus time
for an overall integration time of \(T=10^4\pi \).
The analytical solution was computed using the 1st, 3rd and 5th order normal form.
The numerical solution was obtained using a 4th order Runge--Kutta integration scheme with fixed step
size $\delta_t=10^{-2}$. As expected, the difference between the numerical and analytical
solutions decreases as the order of the normal form increases.

\subsection{Exponential stability estimates}\label{sec:estimates}
In this section we present an application of the Theorem of Section 3 to the
sample provided by the differential system \equ{e19}. We first discuss the
smallness conditions required for the parameters (Section 4.3.1) and then we
compute the stability estimates (Section 4.3.2).

\subsubsection{Bounds on the parameters}
The bounds on the parameters $\varepsilon$ and $\mu$ are due to the smallness conditions imposed by the requirements to invert from original to intermediate variables,
to invert from intermediate to new variables, to satisfy the non--resonance condition in the intermediate variables and the non--resonance condition in the new variables.
With reference to the Appendix~A, assuming \(y_0 =\frac{1}{2}\left(\sqrt{5}+1\right)\), $x_0=0$,
$r_0=0.118$, $s_0=0.1$, $\delta_0=0.05$, $K=20$, one finds $\tilde{r}_0=0.113$, $\tilde{s_0}=0.05$, $\tilde{r}_0'=0.056$, $R_0=0.057$, $S_0=0.025$, $a=0.09$
(the parameters are chosen so to optimize the result). Condition \equ{33ter} requires that $\varepsilon \le 1.2\cdot10^{-4}$, while condition \equ{C6} imposes that $\varepsilon \le 7.2\cdot10^{-4}$ (no requirements are needed on $\mu$). Condition \equ{C7} is satisfied provided $\varepsilon \le 1.2\cdot10^{-4}$, $|\mu| \le 2.0\cdot10^{-4}$, while condition \equ{C8} requires that $\varepsilon \le 3.0\cdot10^{-3}$, $\eta \le 4.75\cdot10^{-3}$. In conclusion we obtain that all conditions are satisfied provided that $\varepsilon \le 1.2\cdot10^{-4}$ and $|\mu| \le 2.0\cdot10^{-4}$.

\subsubsection{Stability estimates}

The final step is to implement the estimates derived from the Theorem, keeping in mind that there are no Fourier modes of the form $h_{10}^{>K}$ in the sample \equ{e19}. Let us write the transformation of the second component from original to final variables as $Y=y+T^{(N)}(y,x,t)$. For given $\varepsilon\leq\varepsilon_0$, $|\mu|\leq\mu_0$ and given \(\lambda _0=\max (\varepsilon_0 ,\mu_0 )\), we calculate the constant $C_p$ as $\left\|T^{(N)}\|_{r_0,s_0}\right/\lambda_0$ and we define $r_1=C_p\lambda_0$. Taking $r_2$ of the same order of magnitude as $r_1$, we set $C_t={r_2}/C_Y$ with $C_Y=\left\|G_{N+1}\|_{R_0,S_0}\right/\lambda_0^{N+1}$. For the variation in action space we set $\tilde\rho_0 =2r_1+r_2\lambda_0$. We finally compute $\tau_0$ from \equ{tau0}. In conclusion
we obtain that
\beqno
\left\|y(t)-y(0)\right\|\leq \tilde\rho_0 \quad {\rm for \ any} \quad t\leq T_0\equiv C_te^{K\tau_0} \ .
\eeqno

\begin{table}
\centering
\(
\begin{array}{lllll}
\hline\hline
 N & 2 & 3 & 4 & 5 \\
\hline
\tau_0&0.851&1.277&1.703&2.129\\
 \|G_{N+1}\|_{R_0,S_0} & 1.966\cdot 10^{-7} & 5.147\cdot 10^{-10} & 1.320\cdot 10^{-12} & 2.053\cdot 10^{-15} \\
 \|T_N\|_{r_0,s_0} & 3.815\cdot10^{-4} & 3.819\cdot10^{-4} & 3.819\cdot10^{-4} & 3.820\cdot10^{-4} \\
 C_p & 1.908 & 1.909 & 1.909 & 1.910 \\
 C_Y & 4.915 & 6.433\cdot 10^1 & 8.251\cdot 10^2 & 6.416\cdot 10^3 \\
 C_t & 3.881\cdot10^{-1}& 2.968\cdot 10^{-2} & 2.314\cdot 10^{-3} & 2.977\cdot 10^{-4} \\
 \tilde\rho_0  & 1.145\cdot10^{-3} & 1.146\cdot10^{-3} & 1.146\cdot10^{-3} & 1.146\cdot10^{-3} \\
 T_0 & 9.702\cdot10^6 & 3.710\cdot 10^9 & 1.446\cdot 10^{12} & 9.302\cdot 10^{14} \\
\hline\hline
\end{array}
\)
\caption{Stability results versus the normalization order $N$. For the definition of the constants see the text. The parameters are taken to optimize the stability time: $\varepsilon_0=1.2\cdot10^{-4}$, $\mu_0=2.0\cdot10^{-4}$, $y_0={1\over 2}(\sqrt{5}+1)$, $r_0=0.118$, $s_0=0.1$, $K=20$. Notice that in this table we fix $K$ and we let $\tau_0$ vary.}
\label{tab5}
\end{table}

\begin{table}
\centering
\(
\begin{array}{lllll}
\hline\hline
N & 2 & 3 & 4 & 5 \\
\hline
 K & 8 & 12 & 17 & 21 \\
 \|G_{N+1}\|_{R_0,S_0} & 1.970\cdot 10^{-7} & 5.169\cdot 10^{-10} & 1.329\cdot10^{-12} & 2.072\cdot 10^{{-15}} \\
 \|T_N\|_{r_0,s_0} & 3.815\cdot10^{-4} & 3.819\cdot10^{-4} & 3.819\cdot10^{-4} & 3.820\cdot10^{-4} \\
 \lambda _0 & 2.0\cdot10^{-4} & 2.0\cdot10^{-4} & 2.0\cdot10^{-4} & 2.0\cdot10^{-4} \\
 C_p & 1.908 & 1.909 & 1.909 & 1.910 \\
 C_Y & 4.924 & 6.461\cdot10^1 & 8.307\cdot10^2 & 6.476\cdot10^3 \\
 C_t & 3.874\cdot10^{-1} & 2.955\cdot 10^{-2} & 2.299\cdot 10^{-3} & 2.949\cdot 10^{-4} \\
 \tilde\rho_0  & 1.145\cdot10^{-3} & 1.146\cdot10^{-3} & 1.146\cdot10^{-3} & 1.146\cdot10^{-3} \\
 T_0 & 3.443\cdot10^{6} & 7.828\cdot10^8 & 1.341\cdot 10^{12} & 5.129\cdot 10^{14} \\
\hline\hline
\end{array}
\)
\caption{Stability results versus the normalization order $N$. For the definition of the constants see the text. The parameters are taken to optimize the stability time: $\varepsilon_0=1.2\cdot10^{-4}$, $\mu_0=2.0\cdot10^{-4}$, $y_0={1\over 2}(\sqrt{5}+1)$, $r_0=0.118$, $s_0=0.1$, $\tau_0=2$. Notice that in this table we fix $\tau_0$ and we let $K$ vary.}
\label{tab5b}
\end{table}

\noindent
The actual values of the parameters, as a function of the normalization order $N$, are summarized in Table \ref{tab5}. The parameters $\eps_0$ and $\mu_0$ are taken from the estimates on the smallness of the parameters (see Section 5.3.1), while $K$ was set equal to $20$.\\

\begin{figure}
\begin{center}
\includegraphics[width=12cm,keepaspectratio]{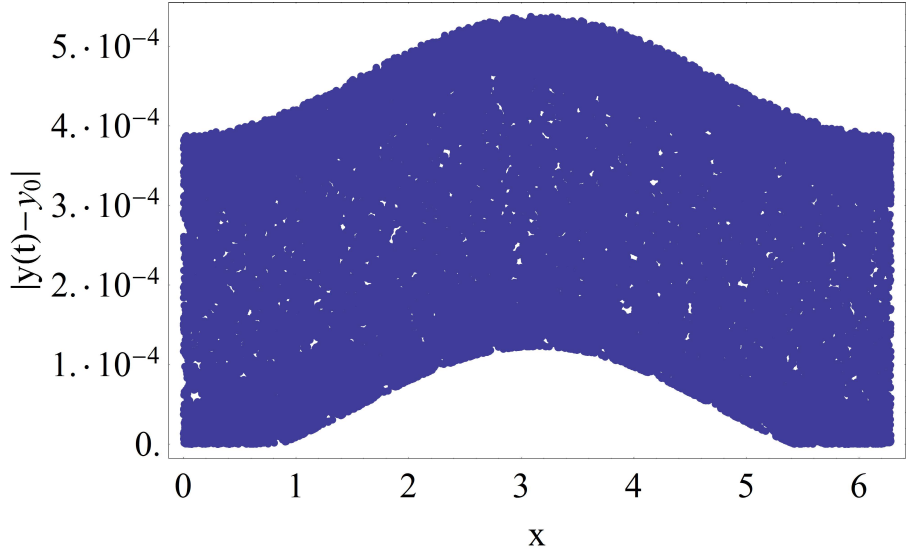}
\caption{Numerical simulation of the original equations of motion  for $3.71\cdot10^9$ integration time. The parameters are taken from Table \ref{tab5}, 3rd column. The numerically found drift in action space is $5.629\cdot10^{-4}$ and it is by a factor $2$ smaller than the upper bound obtained through the analytical estimate, which is equal to $1.146\cdot10^{-3}$.}\label{num}
\end{center}
\end{figure}

\noindent
Similar estimates can be obtained by fixing \(\tau_0\) and calculating $K$ accordingly. Since $K$ also enters into the denominators of the non--resonance conditions \equ{C6} and \equ{C8}, it will also influence the bounds on the smallness of the parameters. The results are shown in Table \ref{tab5b}. Fixing $\tau_0=2$ the choice of $K$ depends on the order of normalization and on the bound on the smallness parameter $\lambda_0$, leading to a slightly different $K \tau_0$ compared to the values given in Table \ref{tab5}. The stability times in Tables \ref{tab5} and \ref{tab5b} are of the same order of magnitude, since the values of $\mu$ and $\varepsilon$ are comparable at all orders.
The stability estimates are checked for the parameters given in Table \ref{tab5b} at order 3 (whose stability time
is compatible with the computer execution time) by comparison with a numerical simulation as shown in Figure \ref{num}. The integration time was set to be of the order of the stability time (i.e., we integrated up to $3.71\cdot10^9$); we found that the deviation of the action is bounded as $5.629\cdot 10^{-4}$, while the analytical estimate provides $1.146\cdot 10^{-3}$ (the numerical deviation is therefore bounded with a safety factor $2$).

\subsection{A system with oscillating energy}\label{sec:oscillating}
\begin{figure}
\begin{center}
\vglue-2cm
\includegraphics[width=7cm,keepaspectratio]{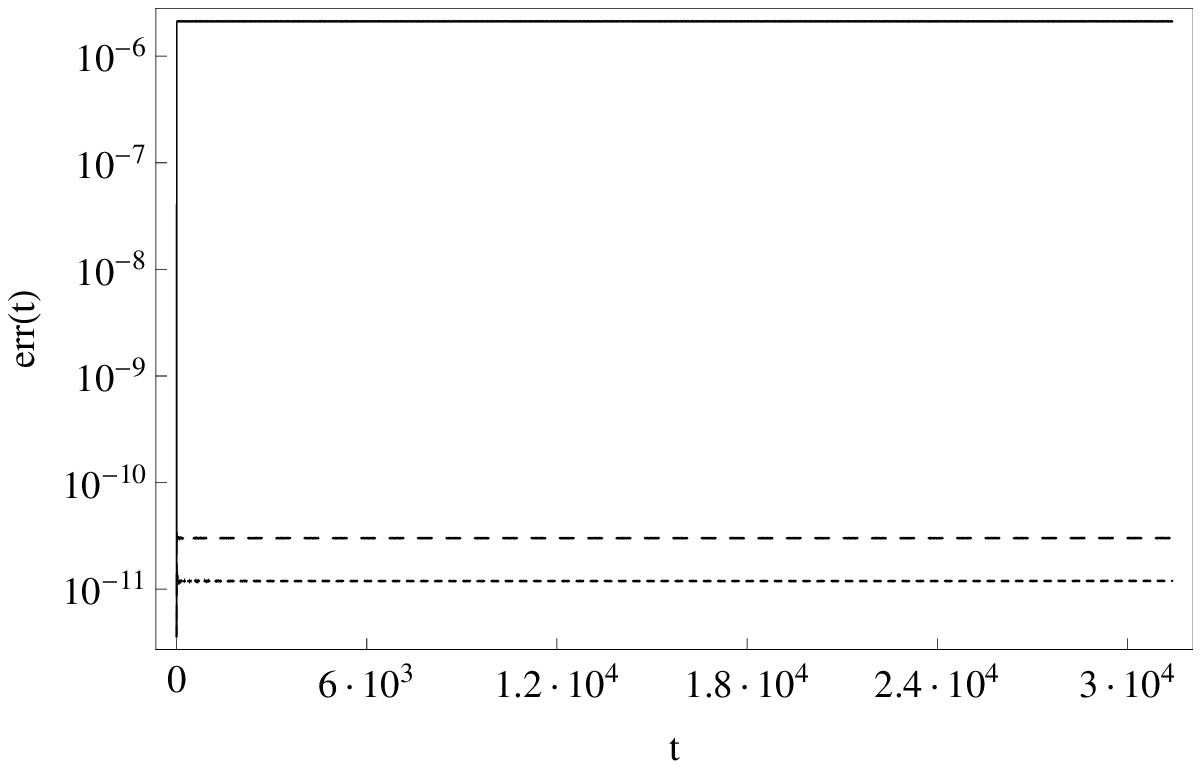}
\hglue-2cm
\includegraphics[width=10cm,keepaspectratio]{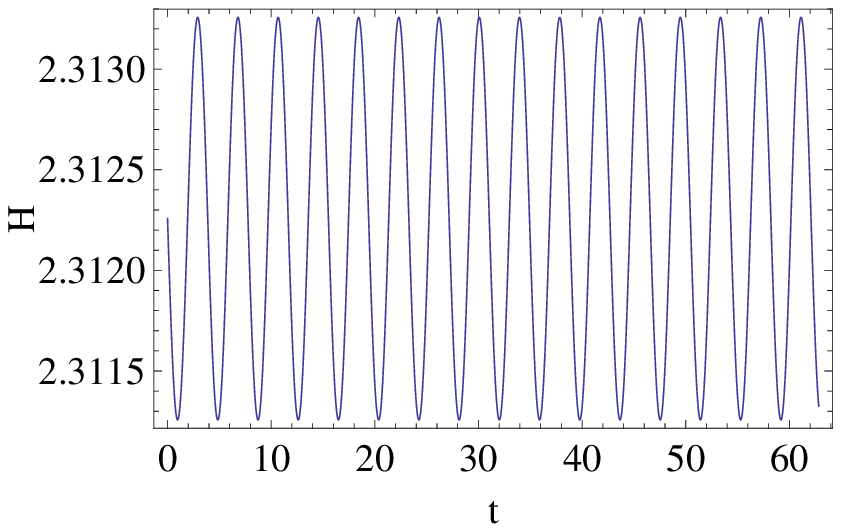}
\caption{Left: relative error between the normal form and the numerical solution. We set $\varepsilon=10^{-3}$, $\mu=10^{-3}$ and $Y_0=\frac{1}{2}\left(\sqrt{5}+1\right)$; the evolution is in good agreement with the solution obtained from the normal form equations (continuous 1st, dashed 3rd, dotted 5th order normal form). Right: behavior of the energy of the system,
which oscillates around a mean value with period $3.86$ .}\label{diss4}
\end{center}
\end{figure}

We conclude by providing an example of a differential system which
admits \sl oscillating energy. \rm To be more precise, we consider the differential equations
\beqano{}
\dot x&=&y \nonumber \\
\dot y&=&-\varepsilon  \left(\sin (x-t)+\sin (x)\right)-\mu  (y \sin (x)-\eta ). \nonumber
\eeqano
The Hamiltonian function for \(\mu =0\) (in the extended phase space) reads as
\beqno{}
H(y,T,x,t)=\frac{y^2}{2}-\varepsilon \left(\cos (x-t)+\cos (x)\right)+T\ , \nonumber
\eeqno
where \(T\) is the conjugated action to the time \(t\). For \(\mu \neq 0\) we get that the variation of the energy is given by
\beqno{}
\frac{d {H}}{d t}=\frac{\partial {H}}{\partial y}\dot y+\frac{\partial {H}}{\partial x}\dot x+\frac{\partial{H}}{\partial T}\dot T+\frac{\partial {H}}{\partial t}\dot t=\mu y^2 \sin (x)- \mu y \eta\ .
\eeqno
Since the normal form equations will provide that \(\eta =0\), we can conclude that the energy is oscillating. The normal form solution to second order provides the following expressions for the transformations:
\beqano{}
\psi_{10}(\tilde{y},\tilde{x},t)&=&\frac{\sin (\tilde{x}-t)}{\tilde{y}-1}+\frac{\sin (\tilde{x})}{\tilde{y}} \nonumber\\
\psi_{20}(\tilde{y},\tilde{x},t)&=&-\frac{\sin (2 \tilde{x}-t)}{2 (\tilde{y}-1) \tilde{y} (2 \tilde{y}-1)}-\frac{\sin (2 \tilde{x}-2t)}{8 (\tilde{y}-1)^3}-\frac{\sin (t)}{2 (\tilde{y}-1)
\tilde{y}}-\frac{\sin (2 \tilde{x})}{8 \tilde{y}^3}\nonumber\\
\beta_{01}(Y,X,t)&=&\cos (X)\nonumber\\
\alpha_{01}(Y,X,t)&=&\frac{\sin (X)}{Y}\nonumber\\
\beta_{02}(Y,X,t)&=&-\frac{\cos (2 X)}{2 Y}\nonumber\\
\beta_{11}(Y,X,t)&=&\frac{Y \cos (2 X-t)}{2 (Y-1)^2 (2 Y-1)}+\frac{(Y-2) \cos (t)}{2 (Y-1)^2}-\frac{\cos (2 X)}{4 Y^2}\nonumber\\
\alpha_{02}(Y,X,t)&=&-\frac{\sin (2 X)}{4 Y^2}\nonumber\\
\alpha_{11}(Y,X,t)&=&-\frac{\sin (2 X-t)}{2 (Y-1) (2 Y-1)^2}+\frac{(Y-3) \sin (t)}{2 (Y-1)^2}-\frac{\sin (2 X)}{8
Y^3} \ .
\eeqano
The corresponding normal form at second order is given by:
\beqano{}
\dot X&=&Y+\frac{1-3Y+3Y^2-2Y^3}{2(Y-1)^3 Y} \varepsilon ^2+\frac{\varepsilon  \mu }
{2 Y^2}+O_2(\varepsilon,\mu)  \nonumber \\
\dot Y&=&O_2(\varepsilon,\mu)\ .
\eeqano
We compare the solution of the normal form equations with the numerical solution by computing the error as in \equ{errt} for $\varepsilon_0=10^{-3}$, $\mu_0=10^{-3}$ and $Y_0=\frac{1}{2}\left(\sqrt{5}+1\right)$, $X_0=0$. Figure~\ref{diss4} (left panel) shows the error for the 1st (continuous line), 3rd (dashed line) and 5th (dotted line) order respectively, up to time $10^4\pi$; the oscillating behavior of the energy is given in the right panel of Figure~\ref{diss4}. The numerical solution was again obtained using a 4th order Runge--Kutta integration scheme with fixed step size $\delta_t=10^{-2}$. We remark that the difference between the numerical and analytical
solutions decreases as the order of the normal form increases. We conclude by mentioning that the bounds on the small parameters as well as the stability estimates can be determined as in Section~\ref{sec:estimates}.

\section{Appendix A}
We discuss the conditions which must be satisfied by the parameters $\varepsilon$, $\mu$, so that the transformation from original to intermediate variables can be inverted, as well as that from intermediate to final variables; moreover, we give conditions on the parameters so that the non--resonance conditions in the intermediate and final variables are satisfied. Such results rely on the following two lemmas which are proven in \cite{GG}.

\vskip.1in

\noindent
\bf Lemma A.1. \sl Let $y_0\in{\R}^\ell$, $(x,t)\in{\T}^{\ell+1}$, $r_0$, $s_0$, $\delta_0$
($\delta_0<s_0$) be
strictly positive parameters and let $g$ be a vector function holomorphic on the
domain $D(y_0,r_0,s_0)\equiv\{(y,x,t)\in{\C}^{2\ell+1}:\ \|y-y_0\|\leq r_0\ ,\
\max_{1\leq j\leq\ell} |\Im(x_j)|\leq s_0\ ,\ |\Im(t)|\leq s_0\}$. Let us consider the equation
\beq{eqb1}
x'=x+g(y,x,t)\ ;
\eeq
if
\beq{b1}
C\|g\|_{r_0,s_0}e^{2s_0}\delta_0^{-1}<1
\eeq
for some positive constant $C$, then \equ{eqb1} can be inverted as
$$
x=x'+G(y,x',t)\ ,
$$
for a suitable function $G$ such that
$$
\|G\|_{r_0,s_0-\delta_0}\leq \|g\|_{r_0,s_0}\ .
$$
\rm

\vskip.1in

\noindent
\bf Lemma A.2. \sl Let $y_0\in{\R}^\ell$, $(x,t)\in{\T}^{\ell+1}$ and let $r_0$, $s_0$, $\tilde r_0$ be
strictly positive parameters with $r_0'<r_0$; let $g$ be a vector function holomorphic on the
domain $D(y_0,r_0,s_0)\equiv\{(y,x,t)\in{\C}^{2\ell+1}:\ \|y-y_0\|\leq r_0\ ,\
\max_{1\leq j\leq\ell} |\Im(x_j)|\leq s_0\ ,\ |\Im(t)|\leq s_0\}$. Let us consider the equation
\beq{eqb2}
y'=y+g(y,x,t)\ ;
\eeq
if
\beq{b2}
C\|g\|_{r_0,s_0}{1\over {r_0-r_0'}}<1
\eeq
for some positive constant $C$, then \equ{eqb2} can be inverted as
$$
y=y'+G(y',x,t)\ ,
$$
for a suitable function $G$ such that
$$
\|G\|_{r_0',s_0}\leq \|g\|_{r_0,s_0}\ .
$$
\rm

\vskip.1in

\noindent
We remark that a careful evaluation of the constant $C$ in \equ{b1} and \equ{b2} shows that it
can be fixed as $C=70$.

\subsection{Inversion of the conservative transformation}
Let us recall the transformation \equ{A} as
\beqa{Ashort}
\tilde x&=&x+\psi_y^{(N)}(\tilde y,x,t)\nonumber\\
y&=&\tilde y+\psi_x^{(N)}(\tilde y,x,t)\ ,
\eeqa
that we wish to invert as
\beqa{Ainv}
x&=&\tilde x+\Gamma^{(x,N)}(\tilde y,\tilde x,t)\nonumber\\
y&=&\tilde y+\Gamma^{(y,N)}(\tilde y,\tilde x,t)\ .
\eeqa
Let $\tilde r_0<r_0$, $\delta_0<s_0$, $\tilde s_0\equiv s_0-\delta_0$; the
inversion of the first in \equ{Ashort} can be performed provided that
$$
70\, \|\psi_y^{(N)}\|_{\tilde r_0,s_0} e^{2s_0}\delta_0^{-1}<1\ ,
$$
with
$$
\|\Gamma^{(x,N)}\|_{\tilde r_0,\tilde s_0}\leq \|\psi_y^{(N)}\|_{\tilde r_0,s_0}\ .
$$
The second in \equ{Ainv} is obtained from
$$
y=\tilde y+\psi_x^{(N)}(\tilde y,\tilde x+\Gamma^{(x,N)}(\tilde y,\tilde x,t),t)\equiv
\tilde y+\Gamma^{(y,N)}(\tilde y,\tilde x,t)\ ,
$$
where
$$
\|\Gamma^{(y,N)}\|_{\tilde r_0,\tilde s_0}\leq \|\psi_x^{(N)}\|_{\tilde r_0,s_0}+
\|\psi_{xx}^{(N)}\|_{\tilde r_0,s_0}\|\Gamma^{(x,N)}\|_{\tilde r_0,\tilde s_0}\ .
$$

\subsection{Non--resonance condition after the conservative normal form}
Taking into account \equ{NR}, we want that the non--resonance condition is satisfied in the
intermediate variables, say for $a>0$:
\beq{omapp}
|\omega(\tilde y)\cdot k+m|>{a\over 2}\ ,\qquad |k|+|m|\leq K\ ,
\eeq
where from \equ{Ashort} we get
\beq{2b}
\tilde y=y+\varepsilon R^{(N)}(y,x,t)\ ,
\eeq
for a suitable function $R^{(N)}$. In fact, the second of \equ{Ashort} can be inverted as in \equ{2b}
provided
$$
70\, \|\psi_x^{(N)}\|_{\tilde r_0,s_0}{1\over {\tilde r_0-\tilde r_0'}}<1\ ,
$$
for $\tilde r_0'<\tilde r_0$ with
$$
\varepsilon \|R^{(N)}\|_{\tilde r_0',s_0}\leq \|\psi_x^{(N)}\|_{\tilde r_0,s_0}\ .
$$
Then we have
\beqno
|\omega(\tilde y)\cdot k+m|
\geq |\omega(y)\cdot k+m|-\varepsilon K\|R^{(N)}\|_{\tilde r_0',s_0} \|\omega_y\|_{r_0}\geq
a-{a\over 2}={a\over 2}\ ,
\eeqno
provided
$$
\varepsilon \leq {a\over {2K\|R^{(N)}\|_{\tilde r_0',s_0}\|\omega_y\|_{r_0}}}\ .
$$

\subsection{Inversion of the dissipative transformation}
Let us now discuss the inversion of \equ{B} that we write for short as
\beqa{3a}
X&=&\tilde x+\alpha^{(N)}(\tilde y,\tilde x,t)\nonumber\\
Y&=&\tilde y+\beta^{(N)}(\tilde y,\tilde x,t)\ ;
\eeqa
we invert \equ{3a} as
\beqano
\tilde x&=&X+\Delta^{(x,N)}(Y,X,t)\nonumber\\
\tilde y&=&Y+\Delta^{(y,N)}(Y,X,t)\ ,
\eeqano
provided $\varepsilon$, $|\mu|$ are sufficiently small. In fact, the first of \equ{3a} can be inverted
provided
$$
70\,\|\alpha^{(N)}\|_{\tilde r_0,\tilde s_0}\ e^{2\tilde s_0}\tilde\delta_0^{-1}<1\ ,
$$
where $\tilde\delta_0<\tilde s_0$. Inverting the equation as
$$
\tilde x=X+ A^{(x,N)}(\tilde y,X,t)\ ,
$$
we have
$$
\|A^{(x,N)}\|_{\tilde r_0,\tilde s_0-\tilde\delta_0}\leq \|\alpha^{(N)}\|_{\tilde r_0,\tilde s_0}\ .
$$
Writing the second of \equ{3a} as
\beqno
Y=\tilde y+\beta^{(N)}(\tilde y,X+A^{(x,N)}(\tilde y,X,t),t)\equiv\tilde y+B^{(y,N)}(\tilde y,X,t)\ ,
\eeqno
we can invert it as
$$
\tilde y=Y+\Delta^{(y,N)}(Y,X,t)\ ,
$$
provided
$$
70\,\|A^{(y,N)}\|_{\tilde r_0,S_0}{1\over {\tilde r_0-R_0}}<1\ ,
$$
with $S_0<\tilde s_0-\tilde\delta_0$, $R_0<\tilde r_0$, being
$$
\|\Delta^{(y,N)}\|_{R_0,S_0}\leq \|A^{(y,N)}\|_{\tilde r_0,S_0}\ .
$$
Notice that $A^{(y,N)}$ can be bounded as
$$
\|A^{(y,N)}\|_{\tilde r_0,S_0}\leq\|\beta^{(N)}\|_{\tilde r_0,\tilde s_0}+
\|\beta^{(N)}_x\|_{\tilde r_0,\tilde s_0}\|A^{(x,N)}\|_{\tilde r_0,S_0}\ .
$$

\subsection{Non--resonance condition after the dissipative normal form}
We now turn to the fulfillment of the non--resonant condition in the new set of variables
$$
|\omega(Y)\cdot k+m|>0\ ,\qquad |k|+|m|\leq K\ .
$$
To this end, we use the transformation
$$
Y=\tilde y+\beta^{(N)}(\tilde y,\tilde x,t;\varepsilon,\mu)
$$
and using \equ{omapp} one easily finds
\beqno
|\omega(Y)\cdot k+m|\geq|\omega(\tilde y)\cdot k+m|-K\|\omega_y\|_{r_0} \|\beta^{(N)}\|_{\tilde r_0,\tilde s_0}
>{a\over 2}-{a\over 4}\ ,
\eeqno
provided the following smallness condition on the parameters is satisfied:
$$
K\,\|\omega_y\|_{r_0} \|\beta^{(N)}\|_{\tilde r_0,\tilde s_0}<{a\over 4}\ .
$$

\vskip.2in

\section{Appendix B}

\noindent
\bf Lemma B.1. \sl Let $f=f(y,x,t)$ be an analytic function on the domain
$C_{r_0}(A)\times C_{s_0}({\T}^{\ell+1})$. Let $f^{>K}(y,x,t)\equiv \sum_{(j,m)\in{\Z}^{\ell+1},
|j|+|m|>K} \hat f_{jm}(y)\ e^{i(j\cdot x+mt)}$ and let $0<\sigma_0<s_0$.
Then, there exists a constant $C_a\equiv C_a(\sigma_0,K)$, such that
\beq{uv}
\|f^{>K}\|_{r_0,s_0}\leq C_a \|f\|_{r_0,s_0+\sigma_0} e^{-(K+1)\sigma_0}\ ,
\eeq
with
\beq{Ca1}
C_a\equiv e^{(K+1){\sigma_0\over 2}}\left(\frac{1+e^{-{\sigma_0\over 2}}}{1-e^{-{\sigma_0\over 2}}}\right)^{\ell+1}\ .
\eeq
\rm

\vskip.1in

\noindent
\bf Proof. \rm From the properties of analytic functions, one has that
$$
|\hat f_{jm}(y)|\leq \|f\|_{r_0,s_0+\sigma_0} e^{-(s_0+\sigma_0)(|j|+|m|)}\ .
$$
Therefore one finds
\beqano
\|f^{>K}\|_{r_0,s_0}&=&\sum_{(j,m)\in{\Z}^{\ell+1}, |j|+|m|>K} |\hat f_{jm}(y)|\
e^{s_0(|j|+|m|)}\nonumber\\
&\leq& \|f\|_{r_0,s_0+\sigma_0}\ \sum_{(j,m)\in{\Z}^{\ell+1}, |j|+|m|>K}
e^{-(s_0+\sigma_0)(|j|+|m|)} e^{s_0(|j|+|m|)}\nonumber\\
&=&\|f\|_{r_0,s_0+\sigma_0}\ \sum_{(j,m)\in{\Z}^{\ell+1}, |j|+|m|>K}
e^{-\sigma_0(|j|+|m|)}\ .
\eeqano
Taking into account that
\beqano
\sum_{(j,m)\in{\Z}^{\ell+1}, |j|+|m|>K} e^{-\sigma_0(|j|+|m|)}&\leq& e^{-(K+1){\sigma_0\over 2}}\
\sum_{(j,m)\in{\Z}^{\ell+1}, |j|+|m|>K} e^{-{\sigma_0\over 2}(|j|+|m|)}\nonumber\\
&\leq& e^{-(K+1){\sigma_0\over 2}}\left(\sum_{p\in{\Z}}e^{-|p|{\sigma_0\over 2}}\right)^{\ell+1}\nonumber\\
&=&e^{-(K+1){\sigma_0\over 2}}\left(\frac{1+e^{-{\sigma_0\over 2}}}{1-e^{-{\sigma_0\over 2}}}\right)^{\ell+1}\ ,
\eeqano
one obtains \equ{uv} with $C_a$ as in \equ{Ca1}. $\Box$

\vskip.2in

\noindent
\bf Acknowledgments. \rm
We deeply thank Luca Biasco, Renato Calleja, Antonio Giorgilli and Jean--Christophe Yoccoz for very useful discussions
and suggestions. We acknowledge the grants ASI
``Studi di Esplorazione del Sistema Solare" and PRIN 2007B3RBEY
``Dynamical Systems and Applications" of MIUR.

\end{document}